\documentclass[12pt]{article}
\title{On groups definable in $p$-adically closed fields}
\author{Anand Pillay\\
\text{ University of Notre Dame} \and
Ningyuan Yao \\
  \text{Fudan University} \and Zhentao Zhang\\
  \text{Fudan University}}

\date{}

\usepackage{amsmath, amssymb, amsthm}    	
\usepackage{fullpage} 	
\usepackage{amscd}
\usepackage{hyperref}
\usepackage[all]{xy}
\usepackage{centernot}
\usepackage{color,xcolor} 
\usepackage{ulem}

\DeclareMathOperator*{\forkindep}{\raise0.2ex\hbox{\ooalign{\hidewidth$\vert$\hidewidth\cr\raise-0.9ex\hbox{$\smile$}}}}

\newcommand{\stab}{\operatorname{stab}}
\newcommand{\Def}{{\operatorname{Def}}}

\newcommand{\GL}{\operatorname{GL}}

\newcommand{\ad}{\operatorname{ad}}
\newcommand{\Ad}{\operatorname{Ad}}

\newcommand{\Aut}{\operatorname{Aut}}

\newcommand{\id}{\operatorname{id}}

\newcommand{\dfg}{\operatorname{dfg}}

\newcommand{\Pow}{\mathcal{P}ow}

\newcommand{\fsg}{\operatorname{fsg}}

\newcommand{\Char}{\operatorname{char}}

\newcommand{\eq}{{\operatorname{eq}}}
\newcommand{\Ker}{\operatorname{Ker}}

\newcommand{\Mat}{{\operatorname{Mat}}}

\newcommand{\gl}{{\mathfrak{gl}}}
\newcommand{\g }{{\mathfrak{g }}}

\newcommand{\C}{\mathbb{C}}
\newcommand{\M}{\mathbb{M}}
\newcommand{\R}{\mathbb{R}}

\newcommand{\Q}{\mathbb{Q}}
\newcommand{\N}{\mathbb{N}}
\newcommand{\Qp}{\mathbb{Q}_p}

\newcommand{\pCF}{p\mathrm{CF}}
\newcommand{\cL}{\mathcal{L}}

\newcommand{\sort}{\mathbf{S}}

\newcommand{\image}{\mathrm{Im}}

\newtheorem{theorem}{Theorem}[section] 
\newtheorem{lemma}[theorem]{Lemma}

\newtheorem{corollary}[theorem]{Corollary}
\newtheorem{fact}[theorem]{Fact}

\newtheorem{conjecture}[theorem]{Conjecture}

\newtheorem{proposition}[theorem]{Proposition}
\newtheorem{proposition-eh}[theorem]{Proposition(?)}
\newtheorem*{theorem-star}{Main Theorem}
\newtheorem*{theorem-star-A}{Theorem A}
\newtheorem*{theorem-star-B}{Theorem B}
\newtheorem*{theorem-star-C}{Theorem C}
\newtheorem*{theorem-star-D}{Theorem D}
\newtheorem*{conjecture-star}{Conjecture}

\newtheorem*{KTC}{Kneser-Tits Conjecture}

\newtheorem*{lemma-star}{Lemma}
\newtheorem*{claim-star}{Claim}

\newtheorem*{claim-star-1}{Claim 1}
\newtheorem*{claim-star-2}{Claim 2}
\newtheorem*{claim-star-3}{Claim 3}
\newtheorem*{claim-star-4}{Claim 4}

\theoremstyle{definition}
\newtheorem{definition}[theorem]{Definition}
\newtheorem{notation}[theorem]{Notation}

\newtheorem{remark}[theorem]{Remark}
\newtheorem{claim}[theorem]{Claim}

\newcommand{\Th}{\text{Th}}

\newcommand{\alg}{\mathrm{alg}}
\newcommand{\Ga}{\mathbb{G}_\mathrm{a}}
\newcommand{\Gm}{\mathbb{G}_\mathrm{m}}
\newcommand{\PSL}{\mathrm{PSL}}
\newcommand{\Lie}{\mathrm{Lie}}
\newcommand{\SL}{\mathrm{SL}}
\newcommand{\sL}{\mathfrak{sl}}

\newcommand{\sq}{\subseteq}

\newenvironment{claimproof}[1][\proofname]
               {
                 \proof[#1]
                 
               }
               {
                 \endproof
               }

\newenvironment{claim1proof}[1][\proofname]
               {
                 \proof[#1]
                 
               }
               {
                 \endproof
               }

\newenvironment{claim2proof}[1][\proofname]
               {
                 \proof[#1]
                 
               }
               {
                 \endproof
               }
\newenvironment{claim3proof}[1][\proofname]
               {
                 \proof[#1]
                 
               }
               {
                 \endproof
               }

\begin{document}
\maketitle

\begin{abstract}
This paper is about the $\dfg$/$\fsg$ decomposition for groups $G$ definable in $p$-adically closed fields. It is proved that for definably amenable $G$, $G$ admits a definable normal $\dfg$ subgroup $H$ such that the quotient $G/H$ is a definable $\fsg$ group. This result was previously known for groups definable in o‑minimal expansions of real closed fields (see \cite{C-P-o-mini}). We also obtain a version for arbitrary (not necessarily definably amenable) groups $G$ definable in $p$-adically closed fields: there exists a definable $\dfg$ subgroup $H$ of $G$ such that the homogeneous space $G/H$ is definable and definably compact. (In the o‑minimal setting this is Fact 3.25 of \cite{Peterzil-Starchenko-mutypes}.) 
Finally, we also isolate the definably amenable part of $G$, which we call the definably amenable component.

Note that $\dfg$ stands for ``has a definable $f$-generic type", and $\fsg$ for ``has finitely satisfiable generics",  which will be discussed together with various equivalences. 

We will need to understand something about groups of the form $G(k)$ where $k$ is a $p$-adically closed field and $G$ a semi-simple algebraic group over $k$, and as part of the analysis we will prove the Kneser-Tits conjecture over $p$-adically closed fields.

\end{abstract}

\section{Introduction}





 We begin with some basic definitions, a summary of earlier work and a statement of our main results. 

First, let $G$ be a group definable in a monster model $\M$ of a complete theory $T$.  Recall from \cite{HPP} that a group $G$ is called \emph{definably amenable} if it admits a  Keisler measure $\mu$ on definable (with parameters) subsets of  $G$ which is invariant under left-translation by elements of  $G$, where a  {\em Keisler measure}  $\mu$  is a finitely additive probability measure on the Boolean algebra  of  definable subsets of $G$.  Definably amenable groups have been studied in detail when the ambient theory $T$ is NIP (does not have the independence property) beginning in \cite{HP-NIP2} and culminating in \cite{Chernikov-Simon}. Two extreme forms of definably amenable groups $G$ have been identified (under an assumption that $T$ is NIP). First, $G$ is said to be $\fsg$ if there is a global complete type $p$ concentrating on $G$ and a small model $M_{0}$ such that every (left) translate of $p$ is finitely satisfiable in $M_{0}$.  Secondly $G$ is said to be (or have) $\dfg$, if there is a global complete type $p$ concentrating on $G$ and a small model $M_{0}$ such that every (left) translate of $p$ is definable over $M_{0}$. In both cases the type $p$  is what is called in \cite{HP-NIP2} a (strongly) $f$-generic type, and the existence of (strongly) $f$-generic types is shown to imply (in fact be equivalent to) the definable amenability of $G$. 

By an $o$-minimal theory we will mean a complete $o$-minimal  expansion of the theory RCF of real closed fields. This is a NIP theory and in fact a special case of a dp-minimal theory ($x=x$ is dp-minimal), see \cite{Simon}.   In \cite{C-P-o-mini} a decomposition theorem was proved for {\em definably connected}  definably amenable group groups definable in an $o$-minimal theory, a combination of Propositions 2.6 and 4.6 there:  $G$ has a normal definably connected definable torsion-free subgroup $H$ such that $G/H$ is definably compact. (In fact $H$ is the greatest definably connected torsion-free subgroup of $G$.) From \cite{HPP} it is known that definable compactness is equivalent to $\fsg$ for definable groups in $o$-minimal theories. On the other hand Proposition 4.7 of \cite{C-P-o-mini} says that definably connected definable torsion-free groups are $\dfg$. So putting it together gives the $\dfg$/$\fsg$ decomposition for definably connected definably amenable groups  (in $o$-minimal theories).  A couple of things follow: first that for definably connected groups, being torsion-free is equivalent to being $\dfg$. And secondly (by passing  from $G$ to $G^{0}$) that the $\dfg$/$\fsg$ decomposition is valid for arbitrary (not necessarily definably connected) groups definably in an $o$-minimal theory.  
Question 1.19 in \cite {PY-question} asks whether the  $\dfg$/$\fsg$ decomposition holds for a definably amenable group in any distal theory. 

Let us now discuss the $p$-adic case. 
Let $p$CF denote the theory of the $p$-adic field $\Qp$ in the ring language. $p$CF is also an example of a distal theory (as is any $o$-minimal theory), and one of the main points of the current paper is to answer positively this Question 1.19 from \cite{PY-question} for the theory $p$CF. (The methods may be able to be applied to so-called P-minimal theories, but we will not do it here). As $p$CF does not have elimination of imaginaries we will distinguish between definable and interpretable. 
Our proof of the $\dfg$/$\fsg$ decomposition in $p$CF will make crucial use of the structure of semi-simple algebraic groups over $p$-adically closed fields (see our Theorem B below and the discussion following it). 

A discussion of definable compactness (as well as a definition) appears in \cite{J-fsg} and will be repeated below.  It is pointed out that for definable subsets $X$ of $M^{n}$ (where $M$ is a model of $p$CF, i.e. a $p$-adically closed field), definable compactness is equivalent to being closed and bounded (with  respect to the usual topology where basic opens are of the form $v(x-a) >\gamma$).  On the other hand for a definable group $G$ in $M$, $G$ can be equipped with a definable manifold structure with respect to $M$ such that  multiplication is continuous. However $G$ is in general far from being (definably) connected (in contradistinction to the $o$-minimal case). 
In \cite{J-fsg} it is shown that for definable groups $G$ in $p$-adically closed fields, $\fsg$ is equivalent to definable compactness, building on the case when $M = \Qp$ (\cite{OP}).

When $M = \Qp$, (or rather when $G$ is defined over $\Qp$) there is a characterization of definable $\dfg$ groups $G$ (see Fact 6.9 in \cite{JY-1}, but proved in \cite{PY-dfg}): $G$ has a subnormal sequence $G_0  \lhd G_1\lhd\cdots \lhd G_n$  where $G_{n}$ is finite index (and normal) in $G$, $G_{0}$ is finite and each $G_{i+1}/G_{i}$ is definably isomorphic to the additive group or a finite index subgroup of the multiplicative group. 
Moreover working over arbitrary models $M$ of $p$CF, it is pointed out in \cite{JY-1} that any one-dimensional definable group which is not definably compact is $\dfg$. 

The $\dfg$/$\fsg$ decomposition (i.e. a positive answer to Question 1.19 from \cite{PY-question}) was proved in \cite{JY-2} when $G$ is abelian, using (among other things) a $p$-adically closed field version of the Peterzil-Steinhorn theorem, proved for abelian groups in \cite{JY-2} (and later extended to arbitrary groups in \cite{JY-3}).

In this paper, we generalize  the $\dfg$/$\fsg$ decomposition from commutative groups to any definably amenable groups:

\begin{theorem-star-A}[=Theorem \ref{theorem-DA-has-dfg-fsg-decomp}]\label{main}
Let $G$ be a definable group in a $p$-adically closed field $M$ which is definably amenable. Then $G$  has a $\dfg$/$\fsg$ decomposition: There is a $M$-definable short exact sequence of $M$-definable groups
$$1 \rightarrow H \rightarrow G \rightarrow C \rightarrow 1$$
where  $H$ is a normal dfg subgroup of $G $, and $C$ is a definable fsg group.  
\end{theorem-star-A}

A key ingredient in the proof of Theorem A is the following, which is the solution of the Kneser-Tits conjecture for $p$-adically closed fields. See Sections 2.3 and 3.2 for the notation.
\begin{theorem-star-B}[=Theorem \ref{theorem-Kneser-Tits Conjecture}]
Let $M$ be a $p$-adically closed field and $G$ a connected, simply connected, almost $M$-simple and $M$-isotropic algebraic group over $M$. Then $G(M)^{+} = G(M)$.
\end{theorem-star-B}
 Let us briefly explain how Theorem B is used. Johnson and Yao proved the dfg/fsg decomposition theorem for commutative definable groups (Theorem 1.1 in \cite{JY-2}) This is extended by an inductive procedure to solvable groups in Proposition \ref{proposition-def-solvable}.  It remains to deal with the ``semi-simple" case (no infinite commutative normal definable subgroups). In the current situation (definable groups in $p$-adically closed fields $M$), a semi-simple definable group $S$ can be assumed to be an open definable subgroup of the set $G(M)$ of $M$-points of a semi-simple algebraic group $G$ over $M$. Lemma \ref{lemma-simple+amenable is compact} says that assuming $S$ to be definably amenable implies that $S$ must be definably compact (so $fsg$).  We can assume $G$ to be almost $M$-simple and $M$-isotropic (if $G$ is $M$-anisotropic then already $G(M)$, so $S$, is definably compact). Hence we know that $G(M)$ is not definably amenable.  Now it follows from Theorem B and a result from \cite{G. Prasad} attributed to J. Tits, that either $S$ has finite index in $G(M)$, contradicting definable amenability of $S$, or $S$ is definably compact, so $fsg$.  (See Corollary \ref{coro-G+=G0}.)

 This dependence on the fine structure of the rational points of semi-simple algebraic groups suggests that proving the $dfg$/$fsg$ decomposition in more general model-theoretic contexts such as groups definable in distal theories will be difficult unless it can be somehow tied to the existence and properties of definable fields.

Let us illustrate an application of the $\dfg$/$\fsg$ decomposition. In the paper \cite{Chernikov-Simon}, Chernikov and Simon studied NIP definably amenable groups and reformulated a question that was initially put forward by Newelski in \cite{Newelski}. Specifically, they explored whether weakly generic types (which are equivalent to $f$-generic types according to Theorem 1.12 in \cite{Chernikov-Simon}) are the same as almost periodic types within their particular settings (see Question 3.32 in \cite{Chernikov-Simon}). Through a series of works presented in \cite{YZ-1}  and \cite{YZ-3}, Yao and Zhang provided an answer to this question for definable groups over the $p$-adics which admit $\dfg$/$\fsg$ decompositions.
A direct consequence of Theorem A in this current paper and Theorem B in the reference \cite{YZ-3} is  as follows:
\begin{corollary}
    Suppose that $G$ is a definably amenable group definable over $\Qp$. Then weak generics coincide with almost periodics iff either $G$ has dfg or $G$ has only bounded many global weak generic types.
\end{corollary}

 When considering an arbitrary definable group $G$ definable in a $p$-adically closed field $M$, we define a $\dfg$ subgroup $H$ of $G$ to be a  {\em $\dfg$ component} of $G$ if it maximizes the dimension $\dim(H)$. It should be noted that the $\dfg$ components of $G$ are not necessarily unique, even in the case where $G$ is definably amenable. Using Theorem A, we have 

\begin{theorem-star-C}[=Theorem~\ref{theorem-VCT}, Theorem~\ref{theorem-G/dfg-component=quotient def-cp}]
Let $G$ be a definable group in a $p$-adically closed field $M$ and $H$ a definable subgroup of $G$. Then $H$ is a $\dfg$ component of $G$ iff $G/H$ is definably compact. Moreover, the $\dfg$ components of $G$ are \emph{conjugate-commensurable}: if $H'$ is another $\dfg$ component of $G$, then there exists $g\in G$ such that $H^g\cap H'$ has finite index in each of $H^g$ and $H'$.
\end{theorem-star-C}

Let us now mention and discuss some related longstanding problems about groups definable in the ``standard model" $\Qp$ and connections to the current paper. These are roughly about the connection between groups definable in the field $\Qp$ and $p$-adic algebraic groups, where (as we mention below) by a $p$-adic algebraic group we mean precisely the group $G(\Qp)$ of $\Qp$-rational points of an algebraic group $G$ defined over $\Qp$. On the other hand a group $G$ definable in $\Qp$ can be definably equipped with the structure of a $p$-adic analytic group (see Lemma 3.8 of \cite{Pillay-fileds-definable-in-Qp}). These are called semialgebraic $p$-adic analytic groups, or $p$-adic Nash groups  (in analogy with the case of real Nash groups).  In any case we may also use the expression ($p$-adic) semialgebraic for definable in $\Qp$. 
The first question was first asked informally in \cite{Pillay-an application to-p-adic-semi-group} (motivated by the truth of the analogous statement in the real case), and later asked formally in \cite{OP}. 

\begin{conjecture} Any open subgroup of a $p$-adic algebraic group $G$ is definable (semialgebraic). 
\end{conjecture} 

This was proved in \cite{open-subgroup} when $G$ is commutative (and as a corollary for reductive $G$ using additionally work of Prasad (\cite{G. Prasad}). It is possible that similar methods could apply to the case of $G$ definably amenable, making use of our Theorem A, and reducing to the case of $G$ being $dfg$. 

\vspace{2mm}
\noindent
The second problem states informally that ``$p$-adic Nash groups are algebraic". (In the real case this is false, and the problem of classifying real Nash groups is still open.) 
\begin{conjecture} Let $G$ be a $p$-adic Nash group. Then (after passing to a definable subgroup of finite index and possibly quotienting by a finite normal subgroup), $G$ is Nash embeddable in a $p$-adic algebraic group. 
\end{conjecture}

Conjecture 1.3 was proved in \cite{JY-2} for $G$ abelian, by proving that (in a saturated elementary extension) $G^{00} = G^{0}$ and making use of a result of Montenegro, Onshuus, Simon in \cite{MOS} (generalizing Hrushovski's stabilizer theorem). In fact the statement $G^{00} = G^{0}$ was then proved for arbitrary $p$-adic Nash groups in \cite{JY-3}, from which they could deduce  by similar methods that Conjecture 1.3 holds for $G$ definably amenable. In fact this definably amenable case could also be deduced from Theorem A above, as Theorem A gives another proof that $G^{00} = G^{0}$ when $G$ is definably amenable. (The $fsg$ group $C$ is compact and it is known that $C^{00} = C^{0}$. On the other hand, by the proof of Lemma 1.15 of \cite{PY-question} the $dfg$ group $H$ satisfies $H^{00} = H^{0}$. This suffices.)

Finally, we study the \emph{definably amenable components} of $G$, which are the maximal definably amenable subgroups of $G$ containing a dfg component. It is worth noting that a definably amenable component of $G$ is not necessarily a definably amenable subgroup of maximal dimension in $G$ (see Remark \ref{rmk-dim-DA-Not max}).
\begin{theorem-star-D}[= Theorem \ref{theorem-max-DA-subgroup}]
Let $G$ be a definable group in a $p$-adically closed field $M$. A subgroup $B \leq G$ is a definably amenable component of $G$ if and only if $B = N_G(H^0)$ for some dfg component $H$ of $G$. Any two definably amenable components of $G$ are conjugate, and $G$ is definably amenable if and only if $G = B$.
\end{theorem-star-D}
Recent work of \cite{YZ-NC} by Yao and Zhang shows that definably amenable components are closely related to the ``Ellis group conjectures". Let $G$ be a group definable in a structure $M$, and consider the $G$-flow given by the action of $G(M)$ on $S_G(M)$. Its enveloping semigroup $(E(S_G(M)), *)$ is an Ellis semigroup. Let $\mathcal{M}$ be a minimal left ideal of $E(S_G(M))$ and $u \in \mathcal{M}$ an idempotent, then $u * \mathcal{M}$ is a group, called the \emph{ideal group} or Ellis group (by model-theorists). In \cite{Newelski}, Newelski noted that $G/G^{00}$ is a homomorphic image of the Ellis group and suggested one may have equality under some conditions. This was made more precise the first author of the current paper who conjectured that we have equality when $T$ is $NIP$ and $G$ is definably amenable.  Chernikov and Simon proved this holds in \cite{Chernikov-Simon}. In \cite{YZ-NC}, Yao and Zhang show that for a group $G$ definable over $\Qp$, the Ellis group of $G$ is isomorphic to the Ellis group  of its definably amenable component. As a corollary, for groups $G$ definable over $\Qp$, the Ellis group is equal to $G/G^{00}$ if and only if $G$ is definably amenable.


\section*{Acknowledgements}
The first author was partially supported by NSF grants DMS-2054271 and DMS-2502292.\\
The second author was supported by the National Natural Science Foundation of China (General Program, Grant No. 12471001) and the Ministry of Education of China (Grant No. 22JJD110002). \\
During the revision of this manuscript, the authors used an AI assistant for spell-checking and grammar correction, and to help identify several cross-reference errors introduced by the renumbering of results. The AI was not used to generate mathematical arguments or bibliographic citations. All content was verified by the authors, who take full responsibility for the final version.


\subsection*{Notation}

We let $\mathcal{L}$ be the language of rings and $T=\pCF$, the theory of $\Q_p$ in $\mathcal{L}$. A model of $T$ is called a $p$-adically closed field. We usually work in a monster model ${\M}$ which is saturated and $|{\M}|$ is strongly inaccessible. We say that $A\subset {\M}$ is small, if $|A|<|{\M}|$. We usually use $M$ to denote a small model. We say that a definable object is defined ``in'' $M$, if it is a definable object in the structure $M$. If we do not say in which model a definable object is, it is a global object, namely, in ${\M}$.  We say that a definable object is defined ``over'' $M$, if it is defined with parameters from $M$. We write $a\in M$ or ${\M}$ to mean that $a$ is a tuple in $M$ or ${\M}$. 

We will assume familiarity with basic notions such as type spaces, heirs, coheirs, definable types etc. References are \cite{Pozit-Book} as well as \cite{Simon}.

\subsection*{Outline}

In Section \ref{sec:Preliminaries}, we introduce several concepts to be used in this paper, including definable compactness, definable manifolds, $ p $-adic definable groups, $ p $-adic algebraic groups, and Lie algebras.

In Section \ref{sect-Kneser-Tits conjecture}, we prove Theorem B, the Kneser-Tits conjecture for any $ p $-adically closed field. From this, we conclude that any open subgroup of a definably simple group is either definably compact or has finite index. Furthermore, we show that if a definable group is definably simple and not definably compact, then it is simple as an abstract group.

In Section \ref{sect-amebility}, we show that a semi-simple $ p $-adic algebraic group is definably amenable if and only if it is definably compact.

In Section \ref{section-dfg}, we show that the property of being a dfg group is definable: Let $ \{G_a \mid a \in Y\} $ be a definable family of definable groups; then $ \{a \in Y \mid G_a \text{ is a dfg group}\} $ is definable.

In Section \ref{section-act-dfg-definably}, we study the action of a dfg group $ G $ on a definably compact space $ X $, showing that there exists some $ x_0 \in X $ with a finite $ G $-orbit, this yields the conjugate-commensurability of dfg components.

In Section \ref{section-act-dfg-fsg}, we show that every definable action of a dfg group on an fsg group is essentially trivial.

In Section \ref{section-dfg-fsg}, we prove our main results: Theorem A and Theorem C.

Finally, Section~\ref{section-DAC} proves Theorem D: the definably amenable components of $G$ (maximal definably amenable subgroups containing a dfg component) are conjugate, and coincide with $G$ exactly when $G$ is definably amenable.

\section{Preliminaries}\label{sec:Preliminaries}

\subsection{Definable compactness  }
We recall the treatment of definable compactness in \cite{JY-1}. 

Let $X$ be a definable set in an arbitrary structure. A {\em definable topology} on $X$ is  a topology $\tau$, such that some  uniformly definable family of definable sets  forms a basis for $\tau$.  When all data are defined over a set $A$ we speak of an $A$-definable topological space $(X,\tau)$. We can apply the definitions (including the uniformly definable basis) to any elementary extension $N$ of $M$, obtaining an $A$-definable in $N$ topological space $(X(N),\tau(N))$. 

\begin{definition}\label{definition-def-comp}
    Let $(X,\tau)$ be a  definable topological space. We call $(X,\tau)$   \emph{definably compact} if for any definable downwards directed family ${\mathcal{F}}=\{Y_s|\ s\in S\}$ of non-empty closed sets,  $\bigcap \mathcal{F}\neq \emptyset$. A definable subset $D\sq X$ is definably compact if it is definably compact as a subspace.
\end{definition}

\begin{fact}[\cite{JY-1}, Remark 2.3]\label{fact-definable compactness elementary}
Let $N\succ M$. Suppose $(X,\tau)$ is a definable topological space in a structure $M$. Then $(X,\tau)$ is definably compact iff $(X(N),\tau(N))$ is definably compact. 
\end{fact}

\begin{fact}[\cite{JY-1}, Fact 2.2]\label{fact closed-subset-of-dc-is-dc}
  Suppose $(X,\tau)$ is definably compact, $Y$ is a subset of $X$. Let $\tau_Y$ be the induced topology on $Y$. If $Y$ is definable and closed, then $(Y,\tau_Y)$ is also definably compact.
\end{fact}

\subsection{Definable manifolds and definable groups in $p$-adically closed fields}
We now specialize to the case where $M$ is a $p$-adically closed field as in the Notation section above. We may sometimes write $k$ for $M$ (to be consistent with algebraic notation). 
Let $\Gamma_k$ denote the value group and $v: k\to \Gamma_k\cup\{\infty\}$ the valuation map. Let $B_{(a,\alpha)}=\{x\in k| v(x-a)>\alpha\}$ for $a\in k$ and $\alpha\in \Gamma_k$. Then $k$ is a topological field with basis given by the definably family ${\mathcal{B}}=\{B_{(a,\alpha)}|\ a\in k, \alpha\in\Gamma_k\}$. We give $k^{n}$ the product topology and note that any definable  subset of $k^{n}$ (with the induced topology)  is a definable topological space in the sense above.

\begin{fact}[\cite{JY-1}, Lemma 2.4]\label{fact-def-cpt=closed-bounded}
 Let $X$ be a definable subset of $k^n$. Then $X$ is definably compact iff $X$ is closed and bounded.
\end{fact}

Recall that if $X\sq k^n$ is definable, then its {\em  topological dimension}   $\dim(X)$ is the greatest $l \leq  n$ such that some projection from  $X$ to $k^l$ contains an open set. 

These notions extend to ``definable manifolds" over $k$:
An $n$-dimensional definable $C^m$-manifold is a topological space $(X,\tau)$ with a covering by finitely many open sets $V_i$ each homeomorphic to an open definable subset $U_{i}$ of  $k^n$ via $f_i$ such that the transition maps are definable and $C^m$.  In this way $(X,\tau)$ becomes a definable topological space (with suitable identifications). 

\begin{fact}[Remark 2.11 in \cite{JY-1}]\label{fact-def-cp-definable}
Definable compactness is a definable property in $\pCF$: Let $k$ be a $p$-adically closed field, $A$ a definable set in $k$, and  $\{X_a|\ a\in A\}$   a
definable family of $k$-definable manifolds. Then
\[\{ a\in A |\  X_a  \ \text{is deﬁnably compact}\}\]
is definable. 
\end{fact}

Let $G$ be a group definable in  $k$. As pointed out in \cite{PY-dfg}, either adapting the methods of \cite{Pillay-minimal-groups}(Proposition 2.5) or by transfer from the case $k = \Qp$ (see  \cite{Pillay-fileds-definable-in-Qp} Lemma 3.8 ), one sees that for any $m <\omega$, $G$ can be definably equipped with a definable topology $\tau_m$, which is a definable $C^m$-manifold in $k$, and for which the group structure is  $C^m$.

Let $(X,\tau)$ be a $\emptyset$-definable topological space,  in $\Qp$. If $X$ is compact, then it is definably compact. By Fact \ref{fact-definable compactness elementary}, $(X(k),\tau(k))$ is definably compact for any elementary extension $k$ of $\Qp$. 
 
A special case is when $X$ is projective $n$-space over $\Qp$: 
$X=\mathbb{P}^n(\Qp)$, the $n$-dimensional projective space in $\Qp$. The usual covering of $X$ as a variety by affine charts makes a $X$ a (definable) compact $p$-adic manifold, so definably compact. Hence for any elementary extension $k$ of $\Qp$, $X(k) = \mathbb{P}^{n}(k)$ is definably compact, as therefore is any Zariski-closed subset. 
Namely we have:


\begin{fact}\label{fact-projective-varieties-are-definably-compact}
    Let $X$ be projective variety over a $p$-adically closed field, then $X(k)$ is definably compact. 
\end{fact}

\subsection{Algebraic groups over $p$-adically closed fields} 


We work now with arbitrary fields $k$ of characteristic $0$. 
Our notation for algebraic groups is basically as in  Milne \cite{J.S. Milne-alg-group}, which we also refer to for basic facts: by an algebraic group $G$ over  $k$ (or an algebraic $k$-group)  we mean a group object in the category of algebraic varieties over $k$.  We may and will  identify $G$ with its group of $\Omega$-points for $\Omega$ an algebraically closed field containing $k$.  On the other hand we can consider the group $G(k)$ of $k$-rational points of $G$. $G(k)$ will in particular be a group (quantifier-free) definable in the field $(k,+,\times)$.  We may sometimes refer to $G(k)$ as a ``$k$-algebraic group"  (generalizing for example notions such as real algebraic group when $k = \R$ and $p$-adic algebraic group when $k = \Q_{p}$). 

For $G$ an algebraic group over $k$ we write $G^{0}$ for the connected component of $G$ as an algebraic group which is also over $k$.  The algebraic group $G$ will be assumed to be connected unless said otherwise.  

\begin{notation}
For the rest of this paper, $\Omega$  will  denote a sufficiently saturated algebraically closed field with characteristic $0$, and  $k$ will always denote  a subfield of $\Omega$ with $|k|<|\Omega|$. And we also assume that our saturated model $\M\succ k$, is a subfield of $\Omega$.  All the algebraic groups we consider in this paper will be algebraic subgroups of $\GL_n(\Omega)$, some $n$.    
\end{notation}

For  solvable algebraic subgroups, we have:
\begin{fact}[Theorem 16.33 of \cite{J.S. Milne-alg-group}]
 Let $G$ be a connected linear algebraic group over $k$.
  Then $G$ is solvable iff $G_{u}$, its  maximal unipotent algebraic subgroup is normal, and $G /G_u$ is an algebraic torus. Moreover $G_{u}$, so also $G/G_{u}$ are connected and over $k$. 
\end{fact}

A connected linear solvable algebraic group $ G $ over a field $ k $ is said to be {\em split} over  $ k$ if $ G $ admits a subnormal series of algebraic subgroups defined over $ k $, whose quotients are $ k $-isomorphic to either the additive group or the multiplicative group. 
Clearly, $ G $ splits over $ k $ if and only if the algebraic torus $ G/G_u $ (where $ G_u $ denotes the unipotent radical of $ G $) splits over $ k $; this is equivalent to $ G/G_u $ being $ k $-isomorphic to a (finite) product of copies of the multiplicative group.

\begin{fact}[Lemma 2.3 of \cite{JY-2}]\label{fact-dfg-exact-seq}
    Suppose that $\M$ is our saturated model of $\pCF$. Let
    \[
    1 \to A \to B \to C \to 1
    \]
    be a definable short exact sequence of definable groups. Then $B$ has dfg iff $A$ and $C$ do.
\end{fact}

 \begin{remark}\label{rmk-fsg-exact-seq}
The same argument as in Lemma 2.3 of \cite{JY-2} shows that the conclusion of Fact~\ref{fact-dfg-exact-seq} also holds for fsg groups: let $\M$ be our saturated model of $\pCF$, and let
\[
1 \to A \to B \to C \to 1
\]
be a definable short exact sequence of definable groups. Then $B$ has fsg iff $A$ and $C$ do.
\end{remark}

Let $G$ be a connected linear solvable algebraic group over $k$ that splits over $k$. If $k$ is a $p$-adically closed field, then both the additive group $\Ga(k)$ and the multiplicative group $\Gm(k)$ have dfg. By induction on $\dim(G)$ and Fact~\ref{fact-dfg-exact-seq}, it follows that $G(k)$ has dfg. We summarize this as follows.
\begin{fact}\label{fact-solvable-split-dfg}
    Let $k$ be a $p$-adically closed field, and let $G$ be a connected linear solvable algebraic group over $k$ that splits over $k$. Then $G(k)$ has dfg.
\end{fact}
In \cite{PY-dfg}, we establish the converse for $k = \Qp$: if $G$ is a connected algebraic group over $\Qp$, then $G(\Qp)$ has dfg if and only if $G$ is a split solvable linear algebraic group over $\Qp$. In this paper, we generalize this equivalence to arbitrary $p$-adically closed fields; see Theorem~\ref{thm-dfg-split}.

Let $G$ be an algebraic group over $k$. 
$G$ is said to be \underline{semi-simple} if the solvable radical $R (G)$ of $G$ is trivial. Here $R (G)$ is the maximal connected normal solvable algebraic subgroup of $G$ (which will be defined over $k$ if $G$ is).  Equivalently  (connected) $G$ is semi-simple if it has no nontrivial connected normal commutative algebraic subgroup. 
We say that $G$ is \underline{$k$-simple} if it is non-commutative and has no proper normal algebraic subgroups over $k$ except for $\{\id_G\}$. It is \underline{almost $k$-simple} if it is  noncommutative and has no
 proper normal infinite connected algebraic subgroup over $k$.  
Note that $G$ can be $k$-simple (resp. almost $k$-simple) without being $k^\alg$-simple (resp. almost $k^\alg$-simple), where $k^\alg$ is the field-theoretic algebraic closure of $k$. 
It is well-known that:
\begin{fact}\label{fact-connected-semi-simples}
Let  $G$ be a connected algebraic group over $k$,  then $G$ is semi-simple iff $G$ is an almost direct product of almost $k$-simple algebraic groups over $k$.
\end{fact}

Recall that an algebraic group $G$ over $k$ is said to be 
{\em reductive} if it is an almost direct product of a semi-simple group and an algebraic torus.
If  $G$ is over $k$ then both the semi-simple part and algebraic torus part are also over $k$.  Suppose $G$ is reductive and over $k$. We call $G$ {\em $k$-isotropic} if it contains a $k$-split algebraic torus over $k$, namely an algebraic subgroup over $k$ isomorphic over $k$ to a finite power of the multiplicative group. {\em $k$-anisotropic} means not $k$-isotropic.

\begin{remark}
Suppose that $k$ is a henselian valued field with $\mathcal{O}_k$ its valuation ring. Let $X \sq k^n$, we say that $X$ is bounded if there is $\alpha\in K$ such that $X\sq \alpha \mathcal{O}_k^n$. We have already seen (Fact \ref{fact-def-cpt=closed-bounded}) that if $X\subseteq k^{n}$ is definable in the $p$-adically closed $k$, then 
$X$ is definably compact iff $X$ is closed and bounded. We view $\GL_{n}(k)$ as a closed subset of $k^{n^2+1}$. 
\end{remark}

\begin{fact}[Theorem 3.1 of \cite{PR-Book}]\label{fact-quotient-by-max-split-subgroup}
Let $H\leq G$ be linear algebraic groups over $\Qp$. 
Then $G(\Qp)/H(\Qp)$ is compact iff $H$ contains a maximal connected $\Qp$-split solvable subgroup of the connected component $G^{0}$ of $G$. In particular $G(\Q_p)$ is compact if and only if $G^{0}$ is reductive and  $\Qp$-anisotropic. 
\end{fact}

\begin{remark}\label{rmk-GL/T-definably compact}
  Let $\mathbb{T}_n(\Omega)$ be the subgroup of $\GL_n(\Omega)$ consists of all upper triangular matrices, then $\mathbb{T}_n(\Omega)$ is a maximal connected $\Qp$-split solvable subgroup of $\GL_n(\Omega)$.   By Fact \ref{fact-quotient-by-max-split-subgroup} we see that $\GL_n(\Qp)/\mathbb{T}_n(\Qp)$ is compact and thus $ \GL_n(k)/\mathbb{T}_n(k)$ is definably compact for any $k\equiv \Qp$.   
\end{remark}

The following Fact is 
the ``Main Theorem" of Prasad's paper \cite{G. Prasad} and is attributed to Bruhat, Tits, and Rouseau. 
\begin{fact}\cite{G. Prasad}\label{fact-semi-simple-unbounded=isotropic}
If $k$ is a Henselian valued field and $G$ is a reductive algebraic group over $k$, then $G(k)$ is bounded if and only if $G$ is $k$-anisotropic.  When $k$ is  $\Q_{p}$ this is a special case of Fact \ref{fact-quotient-by-max-split-subgroup}
\end{fact}

\begin{corollary}\label{coro-definably-compact=not-split}
   Suppose that $k$ is a $p$-adically closed field, and $G$ is a reductive algebraic $k$-group (i.e. algebraic group over $k$), then $G$ is $k$-anisotropic iff $G(k)$ is definably compact.   
\end{corollary}
\begin{proof}
Since $G$ is linear, $G(k)$ is a closed subgroup of some $\GL_n(k)$. So $G(k)$ is definably compact iff $G(k)$ is bounded (Fact \ref{fact-def-cpt=closed-bounded}). We see from Fact \ref{fact-semi-simple-unbounded=isotropic} that $G(k)$ is definably compact iff it is $k$-anisotropic.
\end{proof}




\subsection{Lie algebras}
$k$ remains a field of characteristic $0$.
Recall that a vector space $\g$ over $k$ with a bilinear operation $[ \ ,\  ]$ : $\g\times \g\to \g$ is a \underline{Lie algebra} over $k$ if $[x, x] = 0$ for all $x\in \g$, and  $[ \ ,\  ]$  satisfies the Jacobi identity:
\[
[x, [y, z]] + [y, [z, x]] + [z, [x, y]] = 0
\]
for all $x,y,z\in \g$.  A Lie subalgebra of $\g$ is a subspace $\g_1$ of $\g$ such that $[x, y] \in \g_1$ for all $x,y\in \g_1$. A subspace $I$ of $\g$ is called an ideal of $\g$ if $[x, y] \in I$ for all $x\in I$ and $y\in \g$. Clearly, an  ideal of $\g$ is a subalgebra of $\g$. An ideal $I$ of $\g$ is abelian if $[x, y] =0$ for all $x,y\in I$. We call $\g $  {\em simple} if it has no ideals except itself and $\{0\}$, and call it {\em semi-simple} if it has no abelian ideals except for $\{0\}$.

 \begin{fact}[Corollary 4.15 \cite{J.S.Milne-Lie alg}]\label{fact-semi-simple-model-inv}
    Suppose that   $\g$ is a Lie algebra over $k$, and $k'/k$ is a field extension, then $\g$ is semi-simple iff $\g_{k'}=\g\otimes_kk'$ is semi-simple as a Lie algebra over $k'$.
\end{fact}

 \begin{fact}[Theorem 4.17 ,\cite{J.S.Milne-Lie alg}]\label{fact-semi-simple=prod-of-simples}
   Let   $\g$ be a Lie algebra over $k$, then $\g$ is semi-simple iff it is a direct product of finitely many simple Lie subalgebras $\g_1,...,\g_n$ over $k$, namely, 
    \[
    \g=\g_1\times \cdots \times\g_n,
    \]
    where each $\g_i$ is also a minimal ideal of $\g$.
\end{fact}

As pointed out in Remark 2.2 of \cite{PPS-trans-simple-group}, if $(\g, [\ ,\ ])$ is a Lie algebra over $k$ where $\g$ is a subspace of $k^n$, then $(\g, [\ ,\ ])$ is definable in $k$. Moreover, quantification over subalgebras of $\g$ is first-order, since we actually quantify over bases of these subalgebras. For a subfield $k_0$ of $k$, we say that $\g$ is defined over $k_0$ if there is a basis $v_1, \ldots, v_l \in k_0^n$ of $\g$ such that $[v_i, v_j] \in k_0^n$ for each $i, j \leq l$. Clearly, if $\g$ is defined over $k_0$, then $\g = \g(k_0) \otimes_{k_0} k$. Similarly, if $k'$ is a field extension of $k$, then $\g_{k'} = \g \otimes_k k'$ is $\g(k')$, the realizations of $\g(x)$ in $k'$.

Let $\Mat_n(k)$ be the set of all $n\times n$ matrices over $k$. For each $X,Y\in \Mat_n(k)$, Let $[X,Y]=XY-YX$, then $(\Mat_n(k), [ \ ,\  ])$ is a Lie algebra, denoted by $\gl_n(k)$. 
It is well-known from the classification of the simple Lie algebras over the field $\C$ of complex numbers that:

\begin{fact} For suitable $n$, 
there are finitely many simple Lie subalgebras  $\g_1,...,\g_m$ of $\gl_n(\C)$ such that each $\g_i$ is defined over $\Q$ and every simple Lie subalgebra over $\C$ is isomorphic to some $\g_i$.
\end{fact}

\begin{claim}\label{Claim-finite-many-simple-ones}
    There are finitely many simple Lie subalgebras  $\g_1,...,\g_m$ of $\gl_n( \Omega)$ such that each $\g_i$ is defined over $\Q$ and every simple Lie subalgebra of $\gl_n( \Omega)$ is isomorphic to some $\g_i$.
\end{claim}
\begin{proof}
By definability of simplicity and transfer. 
\end{proof}
 
\begin{claim}\label{claim-finite-semi-simple-lie-alg-over-Q} Fix $ n $. 
   There are finitely many semi-simple Lie subalgebras $ \g_1, \ldots, \g_m $ of $ \gl_n( \Omega) $ such that each $ \g_i $ is defined over $ \Q $ and every semi-simple Lie subalgebra of $ \gl_n( \Omega) $ is isomorphic to some $ \g_i $.
\end{claim}
\begin{proof}
   By Fact \ref{fact-semi-simple=prod-of-simples}, each semi-simple Lie subalgebra $ \g $ of $ \gl_n( \Omega) $ is a product of simple ones. By Claim \ref{Claim-finite-many-simple-ones}, there are only finitely many simple subalgebras of $ \gl_n( \Omega) $, up to isomorphism. Since the dimension of $ \g $ is $ \leq n^2 $, there are only finitely many choices of the simple ideals/subalgebras of $ \g $.
\end{proof}

A {\em geometric structure} (see Section
2 of \cite{HP-groups definbale in local fields})  is an infinite structure $M$ such that  in  any model $N$ of $\Th(M)$, algebraic closure is a pregeometry, and for any formula $\phi(x,y)$ ($x$ a single variable, $y$ a tuple) there is $k$ such that in some (any) model $N$ of $\Th(M)$ if $\phi(x,b)$ has finitely many realizations then has at most $k$ realizations.
We also call $\Th(M)$ a geometric theory. In a model $M$ of a geometric theory given a finite tuple $b$ from $M$ and set $A$ of parameters in $M$, $\dim_{M}(b/A)$ is well-defined as the cardinality of a maximal subtuple of $b$ which is algebraically independent over $A$. And for an $A$-definable subset $X$ of $M^{n}$, $\dim(X)$ can be well-defined as $\max\{\dim_{N}(b/A): b\in X(N), M\prec N\}$, where $N$ is $|M|^+$-saturated.


A special case of a geometric structure is what was called a {\em geometric field} in \cite{HP-groups definbale in local fields}.  Here the language $L$ is the language of rings. We define then (following Remark 2.10 in \cite{HP-groups definbale in local fields}), $k$ to be a geometric field, if $k$ is perfect, model-theoretic algebraic closure coincides with field-theoretic algebraic closure in all models of $\Th(k)$, and we have the uniform bound on finite sets in definable families as above (see \cite{Jo-Ye} for details).

In any  case, we see that in a geometric structure $k$, if $k_{0}$ is a subfield and $b$ a tuple then $\dim_{k}(b/k_{0})$ coincides with $\mathrm{tr.deg}(k_{0}(b)/k_{0})$.  And moreover for $X\subseteq k^{n}$ definable, $\dim(X)$ coincides with the algebro-geometric dimension of the Zariski closure $\tilde X$ of $X$ (in $\Omega^{n}$, where $\Omega$ is an ambient algebraically closed field as mentioned earlier).  As a matter of notation we let $\mathrm{Zdim}(-)$ stand for algebro-geometric dimension (of an algebraic subvariety of $\Omega^{n}$)

    \vspace{2mm}
    \noindent
Note that any model of $\pCF$ is an example of a geometric field, and for any definable set $X$, its  ``geometric dimension'' in the above sense coincides with its “topological dimension”  (see \cite{HP-groups definbale in local fields}, Proposition 2.11 and Fact 2.14).


\begin{claim} \label{claim-vector-space-cap-GL}
Let  $k$ be a geometric field.  Suppose that  $V\sq  \Mat_n(\Omega)$ is a vector space defined over $k$ and $V\cap \GL_n(\Omega)$ is non-empty, then  $\dim(V(k)\cap \GL_n(k))= \mathrm{Zdim}(V)$ (which of course is just the vector space dimension of $V$). 
\end{claim}

\begin{proof}
The main point to see is that $\dim(V(k))$ (as a definable set in $k$) is the same as $\mathrm{Zdim}(V)$, which is because, as $V$ is defined over $k$, $V$ can be ``identified" over $k$ with $k^{m}$ where $m = \dim(V)$.  Hence if $b\in V(k')$ with $k\prec k'$ and $\dim_{k'}(b/k) = \dim(V(k))$ then $b$ is a generic point over $k$ of the variety $V$.  As by assumption $V\cap \GL_{n}(\Omega)$ is nonemptyset it must be a nonempty Zariski open subset of $V$ defined over $k$, where $b\in \GL_{n}(\Omega)$. It follows that again $V$ is the Zariski closure of $V(k) \cap \GL_{n}(k)$ which gives the desired conclusion. 
\end{proof}

\begin{claim}\label{claim-iso-over-K-to-iso-over-k}
    Let   $k $ be a geometric  field, and  $\g_1, \g_2$  Lie subalgebras of $\gl_n(\Omega)$ defined over $k$. If $\g_1$ is isomorphic to $\g_2$, then $\g_1$ is isomorphic to $\g_2$ over $k$, namely, there is a $k$ definable isomorphism $\varphi: \g_1\to \g_2$.
\end{claim}
\begin{proof}
    Let $v_1,...,v_l\in \gl_n(k)$ and $w_1,...,w_l\in \gl_n(k)$ be  bases for $\g_1$ and $\g_2$, respectively. Then with respect to these bases, the isomorphism between $\g_1$ and $\g_2$ is given by a non-singular $l \times l$ matrix $\phi$ over $\Omega$ such that $\phi([v_i,v_j])=[\phi(v_i),\phi(v_j)]$ for $1\leq i,j\leq l$.
    
    Let 
    \[
    V =\{\phi\in \Mat_{l}(\Omega)|\ \phi([v_i,v_j])=[\phi(v_i),\phi(v_j)],\ i,j\leq l\}
    \]
    which note is the set of solutions of a set of linear equations with parameters from $v_1,...,v_l,w_1,...,w_l$ $\in k^{n^2}$. So $V$ is a vector space defined over $k$. Because $V\cap \GL_l(\Omega)$ is non-empty, it follows from Claim \ref{claim-vector-space-cap-GL} that $V(k)\cap \GL_l(k)$ is not empty, which means that there is a $k$-definable isomorphism $\varphi: \g_1\to \g_2$.
\end{proof}

\begin{claim}\label{claim-finit-many-semi-simple-lie-alg}
   Let $k$ be a geometric field. Let  $\g_1,...,\g_m$ be the semi-simple Lie subalgebras of $\gl_n(\Omega)$  defined over $\Q$ given by Claim \ref{claim-finite-semi-simple-lie-alg-over-Q}.
   Then  each semi-simple subalgebra  of $\gl_n(k)$ is isomorphic to some $\g_i(k)$ over $k$. 
\end{claim}
\begin{proof}
     If $\g\leq \gl_n(k)$ is a semi-simple subalgebra, then by Fact \ref{fact-semi-simple-model-inv}, $\g\otimes_k\Omega$ is a semi-simple subalgebra of $\gl_{n}(\Omega)$. By Claim \ref{claim-finite-semi-simple-lie-alg-over-Q}, $\g\otimes_k\Omega$  is isomorphic to some $\g_i$. By Claim \ref{claim-iso-over-K-to-iso-over-k},  we see that $\g\otimes_k\Omega$ is isomorphic to some $\g_i$ over $k$, thus  $\g $ is isomorphic to $\g_i(k)$ over $k$.
\end{proof}
In particular, we have
\begin{claim}\label{claim-finit-many-simple-lie-alg}
   Let $k$ be a geometric field.  
   Then each  $\gl_n(\Omega)$ contains only finitely many $k$-simple Lie subalgebras, up to $k$-isomorphic. 
\end{claim}

If $\g=(k^n,[-,-])$ is a Lie algebra,  the group $\Aut(\g)$ of all automorphisms of $\g$ is an
algebraic subgroup of $\GL_n(k)$. We recall a  well-known fact from \cite{PPS-trans-simple-group}

\begin{fact}[Claim 2.8 of \cite{PPS-trans-simple-group}]\label{fact-dim V=dim Aut V}
    If $k$ is a real closed field, and $\g$ a semi-simple Lie algebra whose underlying vector space is a subspace of $k^n$,  then $\dim(\g)$ = $\dim(\Aut(\g))$.
\end{fact}

 \begin{fact}[Section 4 of Chapter 4 in \cite{OV-book}]\label{fact-complex semi-simple Lie algebra}
   Let $k=\C$ and $\g$ a semi-simple Lie algebra  whose underlying vector space is a subspace of $k^n$. Then
   \begin{itemize}
       \item [(i)] $\Aut(\g)$ is a linear algebraic group, and its Lie algebra $\Lie(\Aut(\g))$ is the algebra of derivations $\mathfrak{Der}(\g)$.
       \item  [(ii)] The adjoint representation 
       \[
       \ad: \g\to \mathfrak{Der}(\g), x\mapsto [x,-]
       \]
       is an isomorphism (see page 203, Problem 3 (Corollary) of \cite{OV-book}).
       In particular, 
       \[
       \dim(\g)=\dim(\mathfrak{Der}(\g))=\dim(\Aut(\g)).
       \]
   \end{itemize}
 \end{fact} 

 \begin{remark}\label{remark-semi-simple-Lie-algebra-over-acf}
   The same proof as Claim 2.8 of \cite{PPS-trans-simple-group} shows that Fact \ref{fact-complex semi-simple Lie algebra} is true whenever $k$ is algebraically closed with characteristic $0$.    
 \end{remark}
The following says that Fact \ref{fact-dim V=dim Aut V} is true for any geometric field $k$ of characteristic $0$.

\begin{claim}\label{claim-dim(g)=dim-Aut(g)}
 If $k$ is a geometric field of characteristic $0$, and $\g$ a semi-simple Lie algebra   whose underlying vector space is a subspace of $k^n$,  then $\dim(\g)$ = $\dim(\Aut(\g))$. 
\end{claim}
\begin{proof}
Since $\g_\Omega=\g\otimes_k\Omega$ is semi-simple,  we have $\mathrm{Zdim}(\Aut(\g_\Omega))=\dim(\g_\Omega)$ by Remark \ref{remark-semi-simple-Lie-algebra-over-acf}. Let $v_1,...,v_l\in k^n$ be a basis for $\g $ and 
    \[
    V=\{\psi\in \Mat(l,\Omega)|\ \psi([v_i,v_j])=[\psi(v_i),\psi(v_j)],\ i,j\leq l\},
    \]
     then  $V$ is a vector space defined over $k$. Clearly, $\Aut(\g_\Omega)= V\cap \GL_l(\Omega)$ and $\Aut(\g)=V(k)\cap \GL_l(k)$. By Claim \ref{claim-vector-space-cap-GL}, we have $\dim(\Aut(\g))= \mathrm{Zdim}(V)= \mathrm{Zdim}(\Aut(\g_\Omega))$. Clearly, $\mathrm{Zdim}(\Aut(\g_\Omega))=\mathrm{Zdim}(\g_\Omega)=\dim(\g)$. We conclude that $\dim(\Aut(\g))=\dim(\g)$ as required.
\end{proof} 


\section{$k$-simple algebraic groups and the Kneser-Tits conjecture}\label{sect-Kneser-Tits conjecture}
For a polynomial $f\in \Omega[x_1,...,x_m]$ and $a\in \Omega^{m}$, by $df_{a}$, we mean the linear map $\sum \frac{\partial f}{\partial x_i}|_a x_{i}$. Let    $G\leq \GL_n(\Omega)$ be  a linear algebraic group. Recall from section 2 of \cite{PPS-JA} that the (Zariski) tangent space $T_{\id_G}(G)$ to $G$ at the identity is defined as the common zero set of linear maps $df_{\id_G}$ in $\Omega$, where $f$ ranges over the polynomials generating the ideal of the variety $G$. The tangent space $T_{\id_G}(G)$ of $G$ at the identity is called the Lie algebra of $G$, and we denote it by $\Lie(G)$. We have $\dim(\Lie(G))=\dim(G)$. Moreover, $\Lie(G)$  is a Lie subalgebra of $\gl_n(\Omega)$,  and the Lie bracket is given by $[X,Y]=XY-YX$.   When $G$ is defined over $k$, then $\Lie(G)$ is defined over $k$ and $\Lie(G)=\Lie(G)(k)\otimes_k \Omega$. (see  Section 2.1.3 of \cite{PR-Book}). 


\subsection{$k$-simple algebraic groups}

\begin{remark}\label{remark-surjective-open-image}
Let $G$ and $H$ be linear algebraic groups over $k\models \pCF$ and $f: G\to H$   a surjective homomorphism over $k$, then $\dim(H(k)) = \dim(f(G(k))$ by an easy computation.  (Also note that for linear connected $G$, $G(k)$ is Zariski-dense in $G$.)  Hence $f(G(k))$  is open in $H(k)$.

\end{remark}

\begin{fact}[Proposition 4.1, \cite{J.S.Milne-Lie alg}]\label{fact-group-semi-simple=lie-alg-semi-simple}
Let   $G$ be a connected linear algebraic groups, then $G$ is semi-simple  iff $\Lie(G)$ is semi-simple. 
\end{fact}

Let $G$ be a connected algebraic group over $k$. Consider the map $\Ad: G\to \Aut(\Lie(G))$ with $\Ad(A)(B)=A^{-1}BA$, which is called the \underline{adjoint representation} of $G$.  Since $\Ker (\Ad)$ is the center of $G$ (see 3.15 of \cite{Borel-book} for example),  when $G$ is centreless,  $\Ker (\Ad)$ is trivial, and thus $\Ad: G\to \Aut(\Lie(G))$ is a $k$-embedding. If $G$ is semi-simple, in addition,  $\Ad: G \to \Aut(\Lie(G ))^0$ is surjective (see page 203, Problem 3 (Corollary) of \cite{OV-book}). Thus, we have that

\begin{fact}\label{fact-adjoint representation-semi-simple algebraic group}
 If $G$ is a centreless semi-simple algebraic group, then $\Ad: G \to \Aut(\Lie(G ))^0$ is an isomorphism.  
\end{fact}

The following Fact is easy to see:
\begin{fact}\label{fact-Lie-alg-iso-imp-aut-lie-iso}
    Let $\g_1,\g_2 =k^n$ be Lie algebra. If $\sigma\in \GL_n(k)$ is an isomorphism from $\g_1$  to $\g_2$, then $\sigma^*: \Aut(\g_1)\to\Aut(\g_2), f\mapsto \sigma f\sigma^{-1}$ is an isomorphism (between the groups of $k$-points of the relevant algebraic groups). 
\end{fact}


\begin{theorem}\label{Thm-semi-simple and K-simple}
 Suppose that $k\models \pCF$ and $G$ is a centreless algebraic group over $k$. 
 \begin{itemize}
     \item If $G$ is semi-simple, then $G$ is $k$-isomorphic to a connected semi-simple algebraic group $H$ over $\Q$. 
     \item If $G$ is $k$-simple, then $G$ is $k$-isomorphic to a connected $\Q$-simple algebraic group $H$ over $\Q$. 
 \end{itemize}   
\end{theorem}
\begin{proof}
    Let $A=\Aut(\Lie(G))^0$. By Fact \ref{fact-adjoint representation-semi-simple algebraic group}, $\Ad: G\to A $ is a $k$-isomorphism. Thus, the restriction of $\Ad$ to $G(k)$ is a $k$-definable isomorphism from $G(k)$ to $A(k)$.
    
    By Fact \ref{fact-group-semi-simple=lie-alg-semi-simple}, $\Lie(G)$ is semi-simple. By Claim \ref{claim-finit-many-semi-simple-lie-alg}, there exists a semi-simple Lie subalgebra $\g$ of $\gl_n(\Omega)$ defined over $\Q$, and a $k$-isomorphism $\sigma: \Lie(G)\to \g$. Let $\Aut(\g)^0$ denote the connected component of $\Aut(\g)$, which is also $\Q$-definable. By Fact \ref{fact-Lie-alg-iso-imp-aut-lie-iso}, $A$ is isomorphic to $\Aut(\g)^0$ via a $k$-definable isomorphism $\sigma^*$. The composition of $\Ad$ and $\sigma^*$ defines a $k$-isomorphism $\varphi: G\to \Aut(\g)^0$, proving the first part. 
    
    The second part is immediate. 
\end{proof}

Recall from \cite{PR-Book} that a connected algebraic group $G$ is said to be \underline{simply connected} if, for any connected algebraic group $H$, any isogeny $f: H \to G$ is an isomorphism, where remember that an  {\em isogeny} is a surjective homomorphism of algebraic groups with finite kernel.

\begin{fact}[ (1994-version) Proposition 2.10 of \cite{PR-Book}]\label{fact-univeral-cover}  
 Assume $G$ to be semi-simple and over $k$. Then there is a simply
connected group $\Tilde{G}$ which is also over $k$ and an isogeny  
$f :\Tilde{G} \to G$ over $k$ (whose kernel is central). ($\Tilde{G}$ is called the universal cover of $G$ and $\Ker(f)$ is the fundamental group of $G$)
\end{fact}

\begin{corollary}\label{coro-simply-connected and almost k-simple}
 Suppose that $k \models \pCF$, and $G$ is a simply connected, almost $k$-simple algebraic group over $k$. Then $G$ is $k$-isomorphic to a simply connected, almost $\Q$-simple algebraic group over $\Q$.
\end{corollary}
\begin{proof}
    Let $G_1 = G/Z(G)$, where $Z(G)$ is the center of $G$. Then $G_1$ is a $k$-simple algebraic group over $k$. By Theorem \ref{Thm-semi-simple and K-simple}, $G_1$ is $k$-isomorphic to a $\Q$-simple algebraic group $H$ over $\Q$. Let $\Tilde{H}$ be the universal cover of $H$; then $\Tilde{H}$ is a simply connected, almost $\Q$-simple algebraic group over $\Q$. Now both $\Tilde{H}$ and $G$ are universal covers of $G_1$, so they are isomorphic over $k$.
\end{proof}

\begin{remark}\label{rmk-semi-simple-is-elementary}
 
  Let $k\models \pCF$ and $k_1\succ k$. Let $G$ be a linear algebraic group over $k$. Then by Theorem \ref{Thm-semi-simple and K-simple}  $G$ is $k$-simple iff $G$ is $k_1$-simple. If $G$ is simply connected, then Corollary \ref{coro-simply-connected and almost k-simple} shows that $G$ is almost $k $-simple iff $G$ is  almost  $k_1$-simple.
 
\end{remark}





Recall that a group $G$ definable in $k$ is called \underline{linear} if there is a definable isomorphism $f: G\to H$ with $H$ a definable subgroup of some $\GL_n(k)$. 

 

\begin{fact}[Lemma 2.11, \cite{JY-3}]\label{fact-antilinear-linear}
 Let $k\models \pCF$. If $G$ is a definable group in $k$, then there is a definable short exact sequence
of groups definable in $k$
\[
1\to A\to G\to H\to 1
\]
where $A$ is commutative-by-finite, and $H$ is  linear group. 
\end{fact}

We  call a group $G$ definable  in $k$,  \emph{definably simple},  if it has no proper nontrivial normal definable subgroup. It is easy to see from Fact \ref{fact-antilinear-linear} that:
\begin{claim}\label{claim-simple-linear}
Let $k\models \pCF$. 
If $H$ is definable in $k$ and definably simple, infinite, and noncommutative, then $H$ is  linear and the Zariski closure $G$ of $H$ is almost $k$-simple.
\end{claim}
\begin{proof}
    By Fact \ref{fact-antilinear-linear}, $H$ is linear. Let $G=G(\Omega)$ be the Zariski closure of $H$ (in suitable $\GL_{n}$). So $G$ is defined over $k$. Note that $G$ is connected, for its connected component $G^0 $ is defined over $k$, and by definable simplicity of $H$, contains $H$. So $G = G^0$. Now suppose that $G$ has some proper normal connected subgroup $N$ defined over $k$. Then $N(k)\cap H$ is a normal $k$-definable subgroup of $H$ so is either trivial or all of $H$. The latter case is not possible as then the Zariski closure of $H$ is in $N$. And in the former case $H$ definably embeds in $G/N$,  so $\dim(H) \leq $ the dimension of $G/N$ as an algebraic group which is $< G$ as an algebraic group, a contradiction.
\end{proof}

Note that $G$ in Claim 3.10 will be semi-simple, so has a finite centre $Z$. Then $H$ definably embeds in $(G/Z)(k)$, in which case and $G/Z$ is $k$-simple.

\begin{proposition}\label{Proposition-definably-compact-not-definably-simple}
Let $ k \models \pCF $. If $ H $ is a definable, definably simple, infinite, and non-commutative group in $ k $, then $ H $ is not definably compact.
\end{proposition}
 
\begin{proof}
By Claim \ref{claim-simple-linear}, there exists a $ k $-simple algebraic group $ G $ such that $ H  $ is an open subgroup of $ G(k) $. By Theorem 
\ref{Thm-semi-simple and K-simple}, $ G $ is $ k $-isomorphic to a $ \Q $-simple algebraic group. We may thus assume $ G $ is already defined over $ \Q $.

Suppose $ H $ is defined by the formula $ \psi(x,b) $ with $ b \in k $. Suppose for a contradiction that $H$ is definably compact. Then the monster model $ \M $ satisfies:
\[
\M \models \left( \psi(\M,b) \leq G(\M) \right) \wedge \left( \psi(\M,b) \text{ is definably compact} \right)
\]

By Theorem 3.13 of \cite{PR-Book}, $ G(\Q_p) $ has finitely many conjugacy classes of maximal compact subgroups, denoted by representatives $ O_1, \dots, O_m $. It follows that:
\[
\Qp\models \forall y\bigg((\psi(\Qp,y) \text{\ is definably compact})\to \exists u\bigvee_{i=1}^m  \left( \psi(\Qp,y)^u \leq O_i \right)\bigg).
\]
Thus we have that 
\[
\M \models \bigvee_{i=1}^m \exists u \left( \psi(\M,b)^u \leq O_i(\M) \right)
\]
Take $g\in G(k)$ and $1\leq i\leq m$ such that $ \psi(\M,b)^g \leq O_i(\M) $.

Now, $ O_i $ admits a uniformly definable family of normal open subgroups: there exist $ L $-formulas $ \theta(x,y) $ and $ \tau(y) $ such that $ \{ \theta(\Qp,y) \mid y \in \tau(\Q_p) \} $ forms a basis of open neighborhoods around the identity $ \id_G $, and each $ \theta(\Q_p, v) $ is a normal subgroup of $ O_1 $ for $ v \in \tau(\Q_p) $. Thus:
\[
\Q_p \models \forall z \forall u \left(( \psi(\Q_p,z)^u \text{ is open in }G(\Qp)) \rightarrow \exists y \left( \theta(\Q_p,y) \subseteq \psi(\Q_p,z)^u \right) \right)
\]

By elementarity, $ \M \models \exists y \left(\tau(y)\wedge (\theta(\M,y) \subseteq \psi(\M,b)^g) \right) $. Let $ c \in \tau(\M) $ be such an element; then $ \theta(\M,c)^{g^{-1}} $ is a non-trivial open normal subgroup of $ \psi(\M,b) $. This contradicts the assumption that $ H $ is definably simple.
\end{proof}

\subsection{The Kneser-Tits Conjecture and Abstract Simplicity}
Let $k$ be an arbitrary field (although we could restrict to characteristic $0$). 
Let   $G$ be an almost $k$-simple algebraic group over $k$. Then $G(k)^+$ is by the definition the subgroup of $G(k)$ generated by the unipotent elements of $G(k)$. Since the set of unipotent elements is normal in $G(k)$, $G(k)^+$ is  a normal subgroup of $G(k)$.
 \begin{fact}[See Fact 6.6 of \cite{PPS-JA}]\label{fact-G+-has no inf nor sbgp}
 Suppose that $G$ is almost $k$-simple and $k$-isotropic. Then 
 \begin{enumerate}
     \item [(i)] $G(k)^+$ has no infinite normal subgroups.
     \item [(ii)] If $G$ is simply connected then $G(k)^+$ contains $T(k)$ for every $k$-split torus $T$ of $G$.
 \end{enumerate}
 \end{fact}





Recall that the Kneser-Tits conjecture is:

\begin{KTC}[See Chapter 7 of \cite{PR-Book}]
Let $G$ be a simply connected, almost $k$-simple, $k$-isotropic algebraic group over $k$. Then $G(k)^+=G(k)$.
\end{KTC}

The {Kneser-Tits Conjecture} is true for $k$ the field $\mathbb{R}$ of real numbers (see Proposition 7.6 of \cite{PR-Book}) or the field $\Qp$ (see Theorem 7.6 of \cite{PR-Book}). It was proved in \cite{PPS-JA} that {Kneser-Tits Conjecture}  is also true for $k$ a real closed field.

\begin{theorem}\label{theorem-Kneser-Tits Conjecture}
    The {Kneser-Tits Conjecture}  is true for every $p$-adically closed field k.  
\end{theorem}
\begin{proof}
Let $G$ be a simply connected,  almost $k$-simple, $k$-isotropic algebraic group over $k$. Our proof is similar to the real closed field case in \cite{PPS-JA} (namely finding ways of using transfer from the standard model), with minor variations.

\begin{claim-star-1}
 We may assume that $k$ is very saturated so that $\Qp\prec k$.       
\end{claim-star-1}
\begin{claim1proof}
    The proof is the same as the proof of the Claim 1 of Proposition 6.8 in \cite{PPS-JA}.
\end{claim1proof}

By Corollary \ref{coro-simply-connected and almost k-simple}, $G$ is $k$-isomorphic to a simply connected, almost $\Qp$-simple algebraic group $H$ over $\Qp$. Now $H$ is $k$-isotropic, so $H(k)$ is unbounded by Fact \ref{fact-semi-simple-unbounded=isotropic}, thus $H(\Qp)$ is also unbounded, and by Fact \ref{fact-semi-simple-unbounded=isotropic} again, $H$ is $\Qp$-isotropic. So we may assume that $G$ is already over $\Qp$ and almost $\Qp$-simple. Let $\mu_G(k)$ be the intersection of all $\Qp$-definable open neighborhoods of $\id_G$ in $k$, and $V_G(k)=\mu_G(k)G(\Qp)$  (Elements of $V_G(k)$ are called the ``finite elements'' of $G(k)$ in the proof of Proposition 6.8 in \cite{PPS-JA}).
\begin{claim-star-2}
   $V_G(k)\leq G(k)^+$. 
\end{claim-star-2}
\begin{claim2proof}
 By the truth of {Kneser-Tits Conjecture} for $\Qp$, $G(\Qp)=G(\Qp)^+$. 
By  Lemma 3 of \cite{G. Prasad}, there are two  unipotent subgroups ${U_1},{U_2}\leq G$ over $\Qp$  such that $G(\Qp)^+$ is generated by ${U_1}(\Qp)$ and ${U_2}(\Qp)$. Let $X = U_{1}(\Qp) \cup U_{2}(\Qp)$.  By Lemma 3.1 of \cite{Pillay-an application to-p-adic-semi-group}, there exist $n\in \mathbb{N}$ such that
$[X]_n=\underbrace{X\cdots X}_{n\ \text{times}}$ contains an open neighborhood $O$ of the identity in $G(\Qp)$. 
Now let  $C\sq G(\Qp)$ be compact. Then $C$ is covered by finitely many translates $c_{i}O$ of $O$ by elements $c_{i}\in C$. As $X$ generates $G(\Qp)$, it follows that $C\subseteq [X]_{n_C}$ for some  $n_C\in \mathbb{N}$. Now let $g\in V_{G}(k)$.  Then $g\in C(k)$, for some definable open compact subset $C$ of $G(\Qp)$.  By what we have observed, and transfer, $g\in [X(k)]_{n_{C}}$. As each element of $X(k)$ is unipotent, $V_{G}(k) \sq G(k)^{+}$, proving Claim 2. 
\end{claim2proof}

\begin{claim-star-3}
  $G(k)^+$ contains all semi-simple elements of $G(k)$.  
\end{claim-star-3}
\begin{claim3proof}
 Let $g\in G(k)$ be semi-simple. Then there is a maximal torus $T$ defined over $k$ and containing $g$. $T$ decomposes (over $k$)  as $T'\cdot T''$ where $T'$ is $k$-split and $T''$ is $k$-isotropic. 
 By Fact \ref{fact-G+-has no inf nor sbgp}, $T'(k)\leq G(k)^+$. On the other hand, $T''(k)$ is definably compact by Corollary \ref{coro-definably-compact=not-split}, so closed and bounded.  

 Now by Theorem 3.13 of \cite{PR-Book},  $G(\Qp)$ has finitely many conjugacy classes of maximal 
compact subgroups. Let these conjugacy classes be represented by $O_{1}$,..,$O_{m}$.   The $O_{i}$'s are definable ($p$-adic semialgebraic) and in fact known to be algebraic. 
Any closed, bounded, definable, subgroup of $G(\Qp)$ is conjugate to a subgroup of some $O_{i}$.   This transfers to $G(k)$, to yield that $T''(k)$ is conjugate to a subgroup of some $O_{i}(k)$.  But $O_{i}(k)$ is contained in $V_{G}(k)$ so by Claim 2 is contained in $G(k)^{+}$. But $G(k)^{+}$ is normal in $G(k)$  (by its definition) whereby $T''(k) \leq G(k)^{+}$.  As $T'(k)\leq G(k)^{+}$, so also $g\in G(k)^{+}$. 
\end{claim3proof}

We now complete the proof of  Theorem  \ref{theorem-Kneser-Tits Conjecture}.  Recall that the Jordan decomposition of a linear algebraic group $G$ writes every element $g$ as a product $g_{s}g_{u}$ where $g_{s}$ is semi-simple and $g_{u}$ unipotent. Moreover by for example Corollary 1 in Section 1.4 of \cite{Borel-book}, if $g\in G(k)$ then both $g_{s}$ and $g_{u}$ are in $G(k)$.  It follows by the definition of $G(k)^{+}$ and Claim 3, that $G(k)^{+} = G(k)$. 
\end{proof}

The following is Theorem 5.8(ii) of \cite{open-subgroup}

\begin{fact}\label{fact-Qp-simple-open-subgp}
Suppose $G$ is a semi-simple algebraic group over $\Qp$ which is  almost $\Qp$-simple, and $\Qp$-isotropic. Then  any open subgroup of $G(\Qp)$ is definable (semialgebraic) and is either compact or of finite index.
\end{fact}

We now generalize Fact \ref{fact-Qp-simple-open-subgp} to any $p$-adically closed field $k$, using Theorem \ref{theorem-Kneser-Tits Conjecture}

\begin{corollary}\label{coro-G+=G0}
Let $k$ be a $p$-adically closed field and $G$ an almost $k$-simple, $k$-isotropic group over $k$. Then
 \begin{enumerate}
     \item [(i)] $G(k)^{+}$ is a definable subgroup of $G(k)$ of finite index,
     \item [(ii)] If $H$ is any open subgroup of $G(k)$, then either $H$ contains $G(k)^{+}$ and is definable of finite index, or $H$ is bounded. In the latter case, if $H$ is also definable then it is definably compact.
 \end{enumerate}
 
\end{corollary}
\begin{proof}
 We may assume $k$ is saturated, so it contains $\Q_p$. By Corollary \ref{coro-simply-connected and almost k-simple} we may assume that $G$ is over $\Q_p$. Hence the universal cover $\pi: \tilde{G} \to G$ is over $\Q_p$. By the truth of Kneser-Tits for $\Q_p$, $\tilde{G}(\Q_p)^+ = \tilde{G}(\Q_p)$. By Corollary 3.20 of \cite{Borel-Tits}, as $\pi$ is an isogeny, $\pi(\tilde{G}(\Q_p))$ has finite index in $G(\Q_p)$. As $k$ is an elementary extension of $\Q_p$, $\pi(\tilde{G}(k))$ has finite index in $G(k)$ (and is visibly definable). By Theorem \ref{theorem-Kneser-Tits Conjecture}, $\tilde{G}(k) = \tilde{G}(k)^+$. By Corollary 6.3 of \cite{Borel-Tits}, $G(k)^+ = \pi(\tilde{G}(k)^+)$, so we conclude (i) of the corollary.

Now for the proof of (ii). Let $H$ be an arbitrary open subgroup of $G(k)$. If $H$ contains $G(k)^+$, then by (i) $H$ is definable and of finite index. Otherwise $H \cap G(k)^+$ is a proper open subgroup of $G(k)^+$. By Theorem (T) from \cite{G. Prasad} (attributed to J. Tits), $H \cap G(k)^+$ is bounded. But $H \cap G(k)^+$ has finite index in $H$, whereby $H$ is bounded too. If $H$ happens to be definable, then as mentioned earlier $H$ will be definably compact. This completes the proof of the Corollary.
\end{proof}

From now on, when we say definably simple we include noncommutative. In \cite{PPS-JA}, the authors proved that:
\begin{fact}\label{fact-definably-simple-groups}
If $G$ is a definably simple group  definable in a real closed field, then either $G$ is definably compact or $G$ is abstractly simple.
\end{fact}

We will now prove the $p$-adically closed field version of Fact \ref{fact-definably-simple-groups}. 

\begin{theorem}\label{thm-definably-simple-groups}
Let $H$ be a definably simple group definable in a $p$-adically closed field $k$. Then $H$ is simple as an abstract group. 
\end{theorem}

\begin{proof}
We may assume that $k$ is saturated. From Fact \ref{fact-antilinear-linear} we see that $H$ is linear, namely definably isomorphic to a subgroup of some $\GL_n(k)$.  So we may assume $H$ to be a subgroup of $\GL_n(k)$. Let $G<\GL (n,\Omega)$ be the Zariski closure of $H$. By Claim \ref{claim-simple-linear}, $G$ is almost $k$-simple. 

As $\dim(H) = \dim(G(k))$, $H$ is an open (so closed) subgroup of $G(k)$.  By Proposition \ref{Proposition-definably-compact-not-definably-simple}, $H$ is not definably compact, so $G$ is $k$-isotropic. We can now use Corollary \ref{coro-G+=G0} (i) and (ii), $H$ contains $G(k)^{+}$. But then $G(k)^{+}$ is definable and of finite index in $H$. Hence by definable simplicity of $H$, $H = G(k)^{+}$. But then by Fact \ref{fact-G+-has no inf nor sbgp} (i), $H$ is abstractly simple.


\end{proof}


\begin{remark}
Let $H$ be a  definable group in  a saturated $p$-adically closed field $k$. Note that if $H$ is definable over $\Qp$ and definably compact, then $\mu_H$, the intersection of all $\Qp$-definable open neighborhood of $\id_H$ in $H$, is an infinite normal subgroup of $H$. So Theorem \ref{thm-definably-simple-groups} is not true when $H$ is definably compact.
\end{remark}


\begin{corollary}
Suppose $H$ is definable in the $p$-adically closed field $k$ and is definably simple. Then any proper definable open subgroup of $H$ is definably compact.
\end{corollary}
\begin{proof}
 This follows from    Corollary \ref{coro-G+=G0} and  the proof of Theorem \ref{thm-definably-simple-groups}. Here are some details. 
 In the proof of Theorem \ref{thm-definably-simple-groups}, $H$ is  $G(k)^{+}$ for some almost $k$-simple $k$-isotropic algebraic group $G$ over $k$, where moreover $G(k)^{+}$ is definable and finite index (so open) in $G(k)$. By Corollary \ref{coro-G+=G0}, part (ii), a proper open subgroup of $H$ is an open subgroup of $G(k)$ which does not contain $G(k)^{+}$ and so must be definably compact. 
\end{proof}

\section{Definable amenability of $p$-adic algebraic groups}\label{sect-amebility}
In this section, $k$ will denote a $p$-adically closed field.

\begin{lemma}\label{PSL2-not-def-amenable}
 $\PSL_2(k)$ is not definably amenable.
\end{lemma}
\begin{proof}
This is like   Remark 5.2 in \cite{HPP}.
\end{proof}

The next lemma is similar to Lemma 4.3 in \cite{C-P-o-mini}.

\begin{lemma}\label{def-amenable-alge-subgroup}
Let $H\leq G$  be algebraic groups. If $G(k)$ is definably amenable  then $H(k)$ is definably amenable. 
\end{lemma}
\begin{proof}
Let $\mu$ be a left $G(k)$-invariant Keisler measure on $G(k)$. Notice that $G(k)/H(k)$ interpretably embeds in $(G/H)(k)$ and the latter is definable (not just interpretable) we can consider $G(k)/H(k)$ as a set definable in $k$. As $\pCF$ has definable Skolem functions (see \cite{Skolem}), there is $S\subset G(k)$ which meets each right coset of $H(k)$ in $G(k)$ in exactly one point. We let $\lambda: Y\mapsto \mu(Y\cdot S)$ for every definable $Y\subset H(k)$. Then  check that $\lambda$ is a left $H(k)$-invariant Keisler measure on $H(k)$. So $H(k)$ is definably amenable.
\end{proof}

The following uses ideas as in  Lemma 4.5 in \cite{C-P-o-mini}.


\begin{lemma}\label{lemma-simple+amenable is compact}
Let $G$ be a semi-simple algebraic group over $k$ and $U$ an open definable subgroup of $G(k)$.  If $U$ is definably amenable then $U$ is definably compact. 
\end{lemma}

\begin{proof}
We know that $G$ is an almost direct product of almost $k$-simple algebraic groups $G_{1},...,G_{n}$ over $k$. 
If $U$ is definably amenable so is each $U\cap G_{i}(k)$, and if each $U\cap G_i(k)$ is definably amenable so is $U$.
Hence we may assume that $G$ is almost $k$-simple. 
If $G$ is $k$-anisotropic then $G(k)$ is definably compact, so $U$ will be too.
Otherwise $G$ is $k$-isotropic.  By Corollary \ref{coro-G+=G0} $G(k)^{+}$ is a definable finite index (open) subgroup of $G(k)$ and $U$ either contains $G(k)^{+}$ or is definably compact. So we may assume that $U$ contains $G(k)^{+}$ (so also has finite index in $G(k)$). We will show that $U$ could not be definably amenable. If it were $G(k)$ would be definably amenable.  Note that $G(k)$ contains nontrivial unipotent elements so its Lie algebra $\Lie(G(k))$ contains nontrivial unipotent elements.  
Then by the Jacobson-Morozov Lemma (see p 87 in \cite{J.S.Milne-Lie alg}), $\Lie(G(k))$ contains an ``$\sL_2$''-triple, in particular a subalgebra $L$ isomorphic to $\sL_{2}$. $L$ will be the Lie algebra of an algebraic subgroup of $G(k)$ isogeneous to $\SL_{2}(k)$.  But then by Lemma \ref{def-amenable-alge-subgroup}, $\SL_{2}(k)$ would be definably amenable, contradicting Lemma \ref{PSL2-not-def-amenable}. This completes the proof. 
\end{proof}

\section{Groups with definable $f$-generics}\label{section-dfg}
Let $T$ be a complete $L$-theory, $\M$   a monster model of $T$.


\subsection{Pro-definability}
We recall the notions of pro-definable sets and pro-definable morphisms from \cite{KY-Pro-definable}: 
Let $ (I, \leq )$ be a small, upward-directed partially ordered set, and let $C$
be a small subset of $\M$. A $C$-definable projective system is denoted as  $\{(Y_i, f_{ij} )| \ i\geq j\in I\}$ and is characterized by the following properties:
\begin{enumerate}
    \item [(1)] each $Y_i$ is a $C$-definable set in $\M$;
    \item [(2)] for every $i \geq j \in I$, the mapping $f_{ij} : Y_i \to Y_j$ is $C$-definable;
    \item [(3)] for all $i \geq j \geq k \in I$, $f_{ii}$ is the identity on $Y_i$ and the composition $f_{ik} = f_{jk} \circ f_{ij}$ holds.
\end{enumerate}
A {\em pro-definable set} $Y$ over $C$ is defined as the projective limit $Y :=  \varprojlim_{i\in I} Y_i$ of a $C$-definable projective system $(Y_i, f_{ij} )$. We say $Y$ is pro-definable if it is pro-definable over $C$ for  some small set of parameters $C$. 

Consider two pro-definable sets $Y :=  \varprojlim_{i\in I} Y_i$ and $Z :=  \varprojlim_{j\in J} Z_j$ over $C$ with associated $C$-definable projective systems $(Y_i, f_{ii'} )$ and $(Z_j, g_{jj'} )$ respectively,  
a {\em pro-definable morphism} $\phi: Y\to Z$ over $C$ is specified by a monotonic function   $d: J\to I$  and a family of $C$-definable functions $\{\phi_{ij} : Y_i \to Z_j | i\geq d(j)\}$ such that, for all $j\geq j' \in J$ and all $i\geq i' \in I$ with $i \geq d(j)$ and $i' \geq d(j')$, the following equality holds:  $\phi_{i'j'}\circ f_{ii'}=g_{jj'}\circ \phi_{ij}$. For each $i_0\in I$ and $j_0\in J$, let $\pi_{i_0}:\prod_{i\in I}Y_i\to Y_{i_0}$ and $\pi_{j_0}:\prod_{j\in J}Z_i\to Z_{j_0}$ be projections, then we have that $\phi_{i_0,j_0}\circ\pi_{i_0}=\pi_{j_0}\circ \phi$ whenever $i_0\geq d(j_0)$.



Suppose that $Y :=  \varprojlim_{i\in I} Y_i$ is pro-definable.  We say that $Y$ is {\em strict pro-definable} if for any finite $I_0\sq I$, $\pi_{I_0}(Y)$ is definable, where  $\pi_{I_0}$ is the projection from $ \prod_{i\in I} Y_i$ onto $ \prod_{i\in I_0}Y_i$.  
\begin{fact}\label{fact-strict-pro-def}\cite{GJ-pro-definability}
    $Y$ is strict pro-definable iff for any definable set $Z$ and pro-definable morphism $f:Y\to Z$, the image $f(Y)$ is definable. 
\end{fact}
A {\em (strict) pro-interpretable set} is a (strict) pro-definable set in $\M^{\eq}$.

\subsection{Uniform  definability  of definable types}

We say that $T$ has  uniform  definability  of definable types (UDDT) if: for any $L$-formula $\psi(x,{y_\psi})$, there is an $L$-formula $d_\psi({y_\psi},z_\psi)$, such that for any  $M\models T$ and any definable type $p(x)\in S_x(M)$, there is $c_{p,\psi}\in M^{|z_\psi|}$ such that for any $b\in M^{|{y_\psi}|}$, 
\[
\psi(x,b)\in p(x)   \iff  M\models  d_\psi(b,c_{p,\psi}) .
\]    
We call $d_\psi({y_\psi},z_\psi)$ the ``definition" of $\psi$ and $c_{p,\psi}$ the corresponding parameter of the definable type $p$. As pointed in the proof of Proposition 4.1 of \cite{KOVACSICS-Ye}, 
we may suppose that the map $\psi\mapsto d_\psi$ factors through
Boolean combinations, namely,  for $L$-formulas $\psi(x,{y_\psi})$ and $\phi(x,{y_\phi})$, we have that $d_{\psi\wedge \phi} $ is $d_\psi \wedge d_\phi $ and $  d_{\neg\psi} $ is $\neg d_\psi $.

If $X\sq \M^{|x|}$ is definable, let $S_X^\Def(\M)$ denote the space of definable types on $X$ and   $\ulcorner X\urcorner$ the code of $X$ (in $\M^\eq$). Working in $\M^{\eq}$, let $\sort_\psi=\{\ulcorner d_\psi(\M^{|{y_\psi}|},c)\urcorner|\ c\in \M^{|z_\psi|}\}$, we see that the map
\[
\mathfrak{d}:S_X^\Def(\M)\to \prod_{\psi\in L}\sort_\psi, \ p\mapsto (\ulcorner d_\psi(\M^{|{y_\psi}|},c_{p,\psi})\urcorner)_{\psi\in L}
\]
is injective. We call $(\ulcorner d_\psi(\M^{|{y_\psi}|},c_{p,\psi})\urcorner)_{\psi\in L}$ a code of $p$, and denote it by $\ulcorner p\urcorner$. Let $X^{\Def}=\{\ulcorner p\urcorner |\ p\in S_X^\Def(\M)\}$.
Kovacsics and Ye showed that when $T$ has UDDT, then $X^\Def$ is pro-interpretable for any definable set $X$ (see Proposition 4.1 of \cite{KOVACSICS-Ye}). The main Theorem of \cite{Beautiful pairs} indicates that $\pCF$ has UDDT and $X^\Def$ is strict pro-interpretable. A recent paper of Guerrero and Johnson showed that

\begin{fact}\label{fact-pro-definability}\cite{GJ-pro-definability}
Suppose that $\M\models \pCF$, then for any definable set $X$,    there is a bijection $f: X^\Def\to Y$ where $Y$ is pro-definable and $f$ is pro-interpretable. It is easy to see from  Fact \ref{fact-strict-pro-def} that $Y$ is strict pro-definable.
 
\end{fact}

\subsection{Dfg groups in $p$-adically closed fields}

From now on until the end of this paper, we work in the context of $\pCF$. Unless explicitly stated otherwise, $k$ denotes a model of $\pCF$, and $\M\succ k$ denotes a monster model of $\pCF$. 


%


\begin{lemma}\label{lemma-dfg-type-definable}
   Let  $\mathcal{G}=\{G_a|\ a\in Y\}$ be   definable family of definable groups in $\M$. Then $\{a\in Y| \ G_a\ \text{has dfg}\}$ is type-definable.
\end{lemma}

\begin{proof}
Suppose that $\mathcal{G}$ is defined by an $L$-formula $G(x,v)$, namely, each $G_a$ is defined by the formula $G(x,a)$. Let $G=G(\mathbb{M}^{|x| },\mathbb{M}^{ |v|})$.  

For each $L$-formula $\psi(x,v,y_\psi)$, let $d_\psi(y_\psi,z_\psi)$ be the definition of $\psi$. Then 
\[
G^\Def=\{(\ulcorner d_\psi (\mathbb{M}^{|y_\psi|},c_{p,\psi})\urcorner)_{\psi\in L}\mid p\in S_G^\Def(\mathbb{M})\}.
\]
Let $\sort_\psi =\{\ulcorner d_\psi (\mathbb{M}^{|y_\psi|},c )\urcorner \mid c\in \mathbb{M}^{|z_\psi|}\} $, then $G^\Def\subseteq \prod_{\psi\in L}\sort_\psi$. Since $G^\Def$ is strictly pro-interpretable, each $\pi_\psi(G^\Def)$ is a definable subset of $\sort_\psi$, where $\pi_\psi: \prod_{\varphi\in L}\sort_\varphi \to \sort_\psi $ is the natural projection.

Identify $G_a$ with $G_a\times\{a\}$, so we may assume $G_a$ is a definable subset of $G$. For any $L$-formula $\varphi(x,y_\varphi)$, let  
\[
\varphi^*(x,y_\varphi,u,v):= G(x,v)\wedge G(u,v)\wedge u\cdot \varphi(x,y_\varphi).
\]
Then for each $p\in S_{G_a}(\mathbb{M})$, $\varphi(x,y_\varphi)\in p(x)$ if and only if $\varphi^*(x,y_\varphi,\id_{G_a},a)\in p$. Assume that $d_{\varphi^*}=d_{\varphi^*}(y_\varphi,u,v,z_{\varphi^*})$.  Then, for $p\in S_{G_a}^\Def(\mathbb{M})$, 
\[
d_{\varphi^*} (\mathbb{M}^{|y_\varphi|},\id_{G_a},a,c_{\varphi^*,p})=d_\varphi(\mathbb{M}^{|y_\varphi|},c_{p,\varphi}),
\]
so 
\[
\ulcorner p\urcorner =(\ulcorner d_{\varphi} (\mathbb{M}^{|y_\varphi|},c_{\varphi ,p})\urcorner)_{\varphi\in L}= (\ulcorner d_{\varphi^*} (\mathbb{M}^{|y_\varphi|},\id_{G_a},a,c_{\varphi^*,p})\urcorner)_{\varphi\in L}
\]
and 
\[
G_a^\Def=\{(\ulcorner d_{\varphi^*} (\mathbb{M}^{|y_\varphi|},\id_{G_a},a, c_{\varphi^*,p})   \urcorner)_{\varphi\in L}\mid p\in S_{G_a}^\Def(\mathbb{M})\}.
\]
Now $G_a$ acts on $G_a^\Def$ via $g\cdot \ulcorner p\urcorner=\ulcorner g\cdot p\urcorner$ for $g\in G_a$ and $p\in S_{G_a}^\Def(\mathbb{M})$. Since 
\[
\ulcorner g\cdot p\urcorner=(\ulcorner d_{\varphi^*} (\mathbb{M}^{|y_\varphi|},\id_{G_a},a, c_{{\varphi^*}, g\cdot p})   \urcorner)_{\varphi\in L}=(\ulcorner d_{\varphi^*} (\mathbb{M}^{|y_\varphi|}, g^{-1} ,a,c_{{\varphi^*},   p})   \urcorner)_{\varphi\in L},
\]
we have 
\[
G_a\cdot \ulcorner p\urcorner=\{(\ulcorner d_{\varphi^*} (\mathbb{M}^{|y_\varphi|},g ,a,c_{{\varphi^*},p})   \urcorner)_{\varphi\in L}\mid g\in G_a \},     
\]
 and 
 \[
 \pi_{\varphi^*}(G_a\cdot \ulcorner p\urcorner)= \{(\ulcorner d_{\varphi^*} (\mathbb{M}^{|y_\varphi|},g ,a,c_{{\varphi^*},   p})   \urcorner)\mid g\in G_a\},
 \]
 which is interpretable over $\{a,c_{{\varphi^*},  p}\}$.

By Fact \ref{fact-pro-definability}, there exist a strictly pro-definable set $X=\varprojlim_{i\in I} X_i$ and a pro-interpretable bijection $\theta: G^{\Def}\to X=\varprojlim X_i$. Then there is a monotone function $\chi: I\to L$ and a family of interpretable maps $\{\theta_{\varphi ,i}: \sort_{\varphi }\to X_i\mid \varphi\geq \chi(i)\}$ such that $\theta_{\varphi ,i}\circ\pi_{\varphi }=\pi_i \circ \theta$.

Now $G_a$ has dfg if and only if there exists $\ulcorner p\urcorner\in G^\Def$ such that $p\vdash G_a$ and $G_a\cdot \ulcorner p\urcorner$ is bounded. Clearly, if $G_a\cdot \ulcorner p\urcorner$ is bounded, then for each $i\in I$, $\pi_i(\theta( G_a\cdot \ulcorner p\urcorner))$ is bounded and hence finite (as $\pi_i(\theta( G_a\cdot \ulcorner p\urcorner))$ is definable). Conversely, if each $\pi_i(\theta( G_a\cdot \ulcorner p\urcorner))$ is finite, then $\theta( G_a\cdot \ulcorner p\urcorner)$ is bounded, so $G_a\cdot \ulcorner p\urcorner$ is bounded. Note that 
\[
\pi_i(\theta( G_a\cdot \ulcorner p\urcorner))=\theta_{\varphi^*,i} (\pi_{\varphi^*} (G_a\cdot \ulcorner p\urcorner))
\]
 for each $\varphi\in L$ and $i\in I$ with $\varphi^* \geq \chi(i)$. Thus, $p$ is a dfg type of $G_a$ if and only if $\theta_{\varphi^*,i} (\pi_{\varphi^*} (G_a\cdot \ulcorner p\urcorner))$ is a finite subset of $X_i$ whenever $\varphi^* \geq \chi(i)$.

For each $\varphi(x,y_\varphi)\in L$ and $i\in I$ with $\varphi^* \geq \chi(i)$, 
\[
\bigg\{\theta_{\varphi^*,i}(\{\ulcorner d_{\varphi^*} (\mathbb{M}^{|y_\varphi|},g ,a, b_{\varphi^*} )   \urcorner\mid g\in G_a\}) \mid a\in Y,b_{\varphi^*}\in\mathbb{M}^{|z_{\varphi^*}|}\bigg\}
\]
is a definable family of definable sets. Since $\pCF$ is geometric, there exists an $L$-formula $\eta_{\varphi^*,i}(v,z_{\varphi^*})$ such that $\theta_{\varphi^*,i}(\{\ulcorner d_{\varphi^*} (\mathbb{M}^{|y_\varphi|}, g ,a,b_{\varphi^*} )   \urcorner\mid g\in G_a\})$ is finite if and only if $\mathbb{M}\models \eta_{\varphi^*,i}(a,b_{\varphi^*})$.

For any finite set $D=\{\varphi_1(x,y_1),\dots, \varphi_m(x,y_m)\}$ of $L$-formulas, let  $\Sigma_D(v,z_{{\varphi_1^*} },...,z_{{\varphi_m^*} })$ denote the following formula:
\[
\forall u  \forall y_1...y_m\bigg(Y(v)\wedge G(u,v)\wedge \bigwedge_{i=1}^m d_{\varphi_i^*} (u ,v, y_{ i} ,z_{{\varphi_i^*} })\to \exists x(\bigwedge_{i=1}^m\varphi_i^*(x,u,v,y_{ i})\wedge G(x,v))\bigg)
\]
and  
\[
\Sigma(v,z_{\varphi^*})_{\varphi\in L}=\{\Sigma_D\mid D \text{ is a finite set of } L\text{-formulas}\}.
\]
We see that for any $(a,b_{\varphi^*})_{\varphi\in L}\models \Sigma$, $(\ulcorner d_{\varphi^*}(\mathbb{M}^{|y_\varphi|},\id_G,a,b_{\varphi^*})\urcorner )_{\varphi\in L}\in G_a^\Def$. Let 
\[
\Delta(v,z_{\varphi^*})_{\varphi\in L}=\{Y(v)\wedge\eta_{\varphi^*,i}(v,z_{ \varphi^*})\mid \varphi\in L, i\in I\},
\]
then for any $(a,b_{\varphi^*})_{\varphi\in L}\models\Sigma\cup \Delta$, $(\ulcorner d_{\varphi^*}(\mathbb{M}^{|y_\varphi|},\id_G,a,b_{\varphi^*})\urcorner )_{\varphi\in L}$ is a code of a dfg type of $G_a$. Conversely, if $p$ is a dfg type of $G_a$, then $(a,c_{\varphi^*,p})_{\varphi\in L}\models \Delta\cup\Sigma$. Thus, $\{a\in Y\mid G_a \text{ has dfg}\}$ is defined by the partial type  
\[
\Phi(v)=\exists (z_{\varphi^*})_{\varphi\in L}(\Delta(v,z_{ \varphi^*})_{\varphi\in L}\cup \Sigma(v,z_{ \varphi^*})_{\varphi\in L}).
\]
\end{proof}




Simon \cite{distal} has isolated a notion, distality, meant to express the property that a NIP theory  has “no stable part”, or is “purely unstable”. Examples include any o-minimal theory and
$\pCF$ \cite{gsm,distal}. Note that the distality of a theory implies that no non-algebraic global types can be both definable and finitely satisfiable (in a small model), so there is no infinite definable group which has both dfg and fsg. It is easy to see that: if  $G$ and $H$ are definable dfg and fsg groups in a $p$-adically closed field, respectively, and $f:G\to H$ is a definable homomorphism, then $f(G)$ is finite. 

Below we characterize dfg groups in $p$-adically closed fields. To carry out induction on dimension, it is necessary to work with quotient spaces, and the following fact will be crucial.

\begin{fact}\label{fact-open-dfg-subgroup} 
Let $H$ be a dfg group acting on a definable set $X$. Then the quotient space $X/H$ is definable (see \cite[Theorem 1.12]{GJ-pro-definability}). 

In particular, if $G$ is a definable group and $H \leq G$ is a dfg subgroup, then the coset space $G/H$ is definable. Furthermore, if $H$ is open in $G$, then $G/H$ is finite, and hence $G$ has dfg (see \cite[Corollary 1.13]{GJ-pro-definability}).
\end{fact}

We first consider the case of linear algebraic groups.

\begin{lemma}\label{lemma-linear-alg-group-k-split-iff-dfg}
Let $ H $ be an algebraic group over $ k $. Then $ H $ is $ k $-split solvable if and only if $ H(k) $ has dfg.
\end{lemma}

\begin{proof}
If $ H $ is $ k $-split solvable, then $ H(k) $ has dfg by Fact~\ref{fact-solvable-split-dfg}. 

Conversely, suppose $ H(k) $ has dfg. We first show that $ H $ is linear. By Chevalley's structure theorem, there is a short exact sequence of algebraic groups over $ k $:
\[
1\to A\to H\overset{\pi}{\to} C\to 1,
\]
where $ A $ is  linear and $ C $ is an abelian variety. Note that $ \pi(H(k)) $ is an open subgroup of $ C(k) $. As an abelian variety is projective, $ C(k) $ is definably compact by Fact~\ref{fact-projective-varieties-are-definably-compact}, hence has fsg. As $ H(k) $ has dfg, $ \pi(H(k)) $ also has dfg by Fact~\ref{fact-dfg-exact-seq}. Since a group with both dfg and fsg must be finite in distal theories, $ \pi(H(k)) $ is finite. It follows that $ \dim(C(k)) = \mathrm{Zdim}(C) = 0 $, so $ C $ is trivial. Thus $ H = A $ (as $H$ is connected) is linear.

We now prove $ H $ is $ k $-split solvable by induction on $ \mathrm{Zdim}(H) $. Let $ R(H) $ denote the solvable radical of $ H $, which is a connected normal $ k $-subgroup. This gives a short exact sequence of algebraic groups over $ k $:
\[
1\to R(H)\to H\overset{\pi}{\to} S\to 1,
\]
where $ S = H/R(H) $ is a semi-simple algebraic group over $ k $. Restricting to $ k $-points yields a short exact sequence of definable groups:
\[
1\to R(H)(k)\to H(k)\overset{\pi}{\to} H(k)/R(H)(k)\to 1.
\]
By Fact~\ref{fact-dfg-exact-seq}, both $ R(H)(k) $ and $ H(k)/R(H)(k) $ are dfg groups. Since $ \pi $ is smooth, $ H(k)/R(H)(k) $ is an open subgroup of $ S(k) $. By Fact~\ref{fact-open-dfg-subgroup}, an open dfg subgroup has finite index in the ambient group, so $ S(k) $ also has dfg.

In particular, $ S(k) $ is definably amenable. By Lemma~\ref{lemma-simple+amenable is compact}, the open subgroup $ S(k) \subseteq S(\M) $ must be definably compact, hence has fsg. Since $ S(k) $ has both dfg and fsg, it is finite, so $ S $ is a finite algebraic group. As $ H $ is connected, $ S $ is trivial, i.e., $ H = R(H) $ is solvable.

It remains to show that the connected solvable group $ H $ is $ k $-split. Let $ H_u $ denote the unipotent radical of $ H $; it is $ k $-defined, normal, and automatically $ k $-split in characteristic 0 (admitting a subnormal series with successive quotients isomorphic to $ \mathbb{G}_a $). By Fact~\ref{fact-solvable-split-dfg}, $ H_u(k) $ is a dfg group. Consider the short exact sequence:
\[
1\to H_u\to H\overset{\pi}{\to} T\to 1,
\]
where $ T = H/H_u $ is a $ k $-torus. By the same argument as above, $ T(k) $ is a dfg group.

By \cite[Corollary 6.11]{JY-1}, $ T(k) $ contains a $1$-dimensional $ k $-definable dfg subgroup $ U $. Let $ \bar{U} $ be its Zariski closure in $ T $; then $ \bar{U} $ is a 1-dimensional $ k $-torus, and $ \bar{U}(k) $ has dfg. A 1-dimensional $ k $-torus is either $ k $-isomorphic to $ \mathbb{G}_m $ (split) or $ k $-anisotropic (definably compact, fsg). Since $ \bar{U}(k) $ has dfg, it cannot be compact, so $ \bar{U} \cong \mathbb{G}_m $ is $ k $-split.

The quotient $ T/\bar{U} $ is again a $ k $-torus, and its $ k $-points form a dfg group. By induction on $ \mathrm{Zdim}(T) $, $ T/\bar{U} $ is $ k $-split. Therefore $ T $ itself is $ k $-split, and hence $ H $, as an extension of a split torus by a split unipotent group, is $ k $-split solvable.
\end{proof}

In general, we have

\begin{theorem}\label{thm-dfg-split}
  Let $G$ be a group definable over $k$. Then $G$ is a dfg group if and only if there exists a $k$-split solvable algebraic group $H = H(\Omega)$ over $k$, a finite-index $k$-definable subgroup $G' \leq G$, and a $k$-definable morphism $f \colon G' \to H(\M)$ such that both $\ker(f)$ and $H(\M)/\operatorname{im}(f)$ are finite.
\end{theorem}

\begin{proof}
The right-to-left direction ($\Leftarrow$) is trivial, so we only verify the converse ($\Rightarrow$).

Let $G$ be a dfg group defined over $k$. It is standard that $G^0 = G^{00}$ for dfg groups; see \cite[Fact 2.3]{PY-dfg}. By \cite[Lemma 2.2]{PY-dfg}, there exists an algebraic group $H$ defined over $k$, a finite-index $k$-definable subgroup $G' \leq G$, and a $k$-definable morphism $f \colon G' \to H(\M)$ such that $\ker(f)$ is finite and $\operatorname{im}(f)$ is open in $H(\M)$. By shrinking $G'$ if necessary, we may assume $H$ is connected. Note that $\operatorname{im}(f)$ is an open dfg subgroup of $H(\M)$; it follows from Fact~\ref{fact-open-dfg-subgroup} that $H(\M)$ has dfg, and hence $H$ is $k$-split solvable by Lemma~\ref{lemma-linear-alg-group-k-split-iff-dfg}.
\end{proof}


\begin{fact}[see  Proposition 16.52 of \cite{J.S. Milne-alg-group}.]Let $F$ be any field with $\Char F=0$, $G$ a linear algebraic group over $F$. Then the following are equivalent:
\begin{enumerate}
    \item [(i)] $G$ is $F$-split solvable algebraic group $G$ over $F$.
\item [(ii)] $G$  is trigonalizable over $F$, namely,  it is isomorphic over $F$ to a group of upper triangular matrices in some $\GL_n(\Omega)$.
\end{enumerate}
\end{fact}

\begin{lemma}\label{lemma-dfg-V-definable}
   Let $\mathcal{G}=\{G_t|\ t\in T\}$  a definable family of definable groups. Then $S=\{t\in T| \ G_t\ \text{has dfg}\}$ is $\vee$-definable.
\end{lemma}
\begin{proof}

Suppose that $\mathcal{G}$ is defined by the formula $G(x,v)$. Let $\mathbb{T}_n(\Omega)$ be the group consists of all upper triangular matrices in $\GL_n(\Omega)$.

Let $\psi(y,u_\psi)$ be a quantifier free $\cL$-formula with $|y|=n^2$, let $\theta_\psi(u_\psi)$ be a quantifier free $\cL$-formula such that
\[
\Omega\models \forall u_\psi(\theta_\psi(u_\psi)\leftrightarrow (\psi(\Omega,u_\psi)\ \text{is a subgroup of } \mathbb{T}_n(\Omega))).
\]
For any $b\in \M$, $\M\models \theta_\psi(b)$ iff $\Omega\models \theta_\psi(b)$  since $\theta_\psi $ is quantifier free. So $\psi(\Omega,b)$ is an algebraic subgroup of $\mathbb{T}_n(\Omega)$ iff  $\M\models \theta_\psi(b)$.

For each $\cL$-formula $\varphi(x ,y;v,u_\psi,z_\varphi)$, let $E_{\varphi, \psi,m }(v,u_\psi,z_\varphi)$ be a formula saying that 
$\varphi(x ,y;\ v,u_\psi,z_\varphi)$ defines a homomorphism, denoted by $f_\varphi$, from $G(\M ,v)$ to   $\psi({\M} ,u_\psi)\wedge\theta_\psi(u_\psi)$ such that  
\[
|(\psi({\M} ,u_\psi)\wedge\theta_\psi(u_\psi))/f_\varphi(G(\M ,v))|\leq m,|\ker(f_\varphi)|\leq m.
\]
  Let $\cL^{\mathrm{qf}}$ be the set of all quantifier-free $\cL$-formulas, then 
    \[
    S=\bigcup_{\varphi\in\cL, \psi \in \cL^{\mathrm{qf}},m\in \N^+ }\exists u_\psi\exists z_\varphi E_{\varphi, \psi,m }( {\M}^{|v|},u_\psi,z_\varphi),
    \]
    which is  $\vee$-definable.    
    
\end{proof}

We conclude from Lemma \ref{lemma-dfg-type-definable} and Lemma \ref{lemma-dfg-V-definable}  that being a dfg group is a definable property:
\begin{theorem}\label{thm-dfg-definable}
Let   $G(x,y)$ be an $\mathcal{L}$-formula. Then the set 
\[
D=\{b\in{\M}^{|y|}:G({\M}^{|x|},b)\text{ is a dfg group}\}
\] is $\emptyset$-definable.
\end{theorem}

The characterization of dfg groups also helps us study intersections of dfg subgroups. Understanding these intersections is helpful for seeing how dfg subgroups are distributed within definable groups. First, we show that the intersection of finitely many dfg subgroups is again a dfg subgroup.

\begin{lemma}\label{lemma-dfg-intersection}
Let $G$ be a definable group, and let $H, K \leq G$ be dfg subgroups. Then $H \cap K$ is dfg.
\end{lemma}

\begin{proof}
We proceed by induction on $n := \dim H$.
If $n = 0$, then $H$ is finite and the claim is trivial.

Now suppose $n > 0$. By Theorem~\ref{thm-dfg-split}, there exists a normal dfg subgroup $H_1 \lhd H$ such that $H/H_1$ is a one-dimensional dfg group. Note that $H_1 \cap K$ is normal in $H \cap K$, and by the induction hypothesis it is a dfg group. We have the definable short exact sequence:
\[
1 \to H_1 \cap K \to H \cap K \to \frac{H \cap K}{H_1 \cap K} \to 1.
\]
If $\dim(H_1 \cap K) = \dim(H \cap K)$, then $H_1 \cap K$ is an open dfg subgroup of $H \cap K$, so $H \cap K$ has dfg by Fact~\ref{fact-open-dfg-subgroup}.

Now suppose $\dim(H_1 \cap K) < \dim(H \cap K)$. Set $L := \frac{H \cap K}{H_1 \cap K}$. Note that $L$ embeds naturally into $H/H_1$ as the subgroup $\{hH_1 : h \in H \cap K\}$, and $\dim(L) = \dim(H/H_1) = 1$.

\begin{claim-star}
$L$ has finite index in $H/H_1$.
\end{claim-star}
\begin{claimproof}
Since $K$ has dfg, the quotient $G/K$ is definable by Fact~\ref{fact-open-dfg-subgroup}. Let
$X := \{hK : h \in H\} \subseteq G/K$. Then $X$ is definable and is naturally identified with $H/(H \cap K)$. The group $H$ acts on $X$ by $(h_1, hK) \mapsto h_1 h K$, so the dfg subgroup $H_1$ also acts on $X$. By Fact~\ref{fact-open-dfg-subgroup} again, the orbit space $Y := X/H_1 = \{H_1 h K : h \in H\}$ is definable.

Since $H_1$ is normal in $H$, we have $H_1 h K = h H_1 K$ for all $h \in H$. Define the map
$\rho \colon H/H_1 \to Y$ by $\rho(hH_1) = h H_1 K\ (= H_1 h K)$. It is straightforward to check that $\rho$ is surjective, and its fibers are exactly the left cosets of $L$ in $H/H_1$. Hence $Y$ is a definable copy of $(H/H_1)/L$, and
\[
\dim(Y) = \dim(H/H_1) - \dim(L) = 1 - 1 = 0.
\]
Thus $Y$ is finite, so $L$ has finite index in $H/H_1$.
\end{claimproof}

Since $H/H_1$ is dfg and $L$ has finite index, $L$ is also a dfg group. Therefore $H \cap K$ is an extension of two dfg groups, and hence has dfg by Fact~\ref{fact-dfg-exact-seq}.
\end{proof}

We now prove that the intersection of any definable family of dfg subgroups is essentially a finite intersection.

\begin{lemma}\label{lemma-intersection-of-a-fimily-of-finite-index-subgroup}
   Let $H$ be a definable group  and $\{H_a|\ a\in Y\}$ a definable family of   groups  of $H$. If each $H_a$ has finite index in $H$, then there are finitely many $a_1,...,a_n\in Y$ such that $\bigcap_{a\in Y }H_a=\bigcap_{i=1 }^n H_{a_i}$.  
\end{lemma}
\begin{proof}
Let $ A =\bigcap_{a\in Y }H_a $. Then $A$ is a definable subgroup of $H $. Since  each $H_a $ has finite index in $H $, we see that $H^0\vdash A $. It follows that $H /A $ is finite, and which implies that there are $a_1,...,a_n\in Y$ such that $|H  /(\bigcap_{i=1}^n H_{a_i} )|=|H /A |$. So $\bigcap_{a\in Y }H_a =A =\bigcap_{i=1 }^n H_{a_i} $ as required.
\end{proof}

 Combining the two lemmas above, we have

\begin{theorem}\label{dfg-inters-in-fact-finite}
Let $G$ be a definable group and $\{H_a:a\in Y\}$ a definable family of dfg subgroups of $G$. Then there are finitely many $a_1,\dots, a_n\in Y$ such that $\bigcap_{a\in Y} H_a=\bigcap_{i=1}^n H_{a_i}$ and the intersection is a dfg subgroup.
\end{theorem}
\begin{proof}
We consider all finite intersections of $\{H_a:a\in Y\}$ and pick a finite subset $\{a_1,\dots,a_m\}\subset Y$ such that $H:=\bigcap_{i=1}^m H_{a_i}$ has minimal dimension. Then by Lemma~\ref{lemma-dfg-intersection}, $H$ is a dfg group; and by Fact~\ref{fact-open-dfg-subgroup}, every $H_a\cap H$ is of finite index in $H$.   
It follows from Lemma~\ref{lemma-intersection-of-a-fimily-of-finite-index-subgroup} that there are finitely many $a_{m+1},\dots, a_n\in Y$ such that 
\[
\bigcap_{a\in Y} H_a=\bigcap_{a\in Y} (H_{a}\cap H)=\bigcap_{i=m+1}^n (H_{a_i}\cap H)=\bigcap_{i=1}^n H_{a_i}.
\]
Again, by Lemma~\ref{lemma-dfg-intersection}, the intersection is a dfg subgroup.
\end{proof}
\begin{corollary}\label{coro-intersection-of-a-family-of-dfg}
    Let $G$ be a dfg group definable over $k$, $Y$ a $k$-definable set, and $\{H_a \mid a\in Y\}$ a definable family. If $\{H_a(k) \mid a\in Y(k)\}$ is a family of dfg subgroups of $G(k)$, then there are finitely many $a_1,\dots,a_n\in Y(k)$ such that $\bigcap_{a\in Y} H_a = \bigcap_{a\in Y(k)} H_a = \bigcap_{i=1}^n H_{a_i}$.
\end{corollary}

\begin{proof}
By Theorem~\ref{thm-dfg-definable}, we see that $\{H_a \mid a\in Y\}$ is a definable family of dfg subgroups of $G$. By Theorem~\ref{dfg-inters-in-fact-finite}, there are $b_1,\dots,b_n\in Y(\M)$ such that $\bigcap_{a\in Y} H_a = \bigcap_{i=1}^n H_{b_i}$, namely
\[
\M \models \exists y_1,\dots,y_n\in Y \left( \bigcap_{a\in Y} H_a = \bigcap_{i=1}^n H_{y_i} \right),
\]
which is a first-order sentence. Thus we can find $a_1,\dots,a_n\in Y(k)$ such that $\bigcap_{a\in Y} H_a = \bigcap_{i=1}^n H_{a_i}$. Since $\bigcap_{a\in Y} H_a \subseteq \bigcap_{a\in Y(k)} H_a \subseteq \bigcap_{i=1}^n H_{a_i}$, the equality holds as required.
\end{proof}

\section{Actions of $\dfg$ groups on definably compact sets}\label{section-act-dfg-definably}

Let $X\sq \M^m$ be a definable subset over $k$. Let $n$ be a natural number; the relation $S_n\sq X^n \times X^n$ is defined as follows: for $a=(a_1,\cdots, a_n), b=(b_1,\cdots, b_n)\in X^n$,
\[
S_n(a,b)\iff \text{ there is a permutation } \sigma \text{ on } \{1,\dots, n\} \text{ such that }  b_i=a_{\sigma(i)} \text{ for each } 1\leq i\leq n.
\]
Clearly, $S_n$ is an equivalence relation on $X^n$, and $\Pow_n(X):=X^n/S_n$ can be regarded as the collection of all $n$-element unordered subsets (i.e., $n$-choices) of $X$. For $a_1,\dots,a_n\in X$ with $a_i=(a_{i1},\dots, a_{im})$ for each $i$, map $a=(a_1,\dots,a_n)\in X^n$ to the tuple of coefficients of the polynomial $\prod_i \left(Y+\sum_j a_{ij}X_j\right)$. Denote this map by $\eta$; then $\eta$ is definable, and for any $a,b\in X^n$, $\eta(a)=\eta(b)$ if and only if $S_n(a,b)$. Thus, we can identify $\Pow_n(X)$ with the image of $\eta$. If $X$ is equipped with a definable topology, the topology on $\Pow_n(X)$ is defined by: a subset $O\subset \Pow_n(X)$ is open if and only if $\eta^{-1}(O)$ is open in $X^n$. This makes $\Pow_n(X)$ a definable topological space for which $\eta$ is continuous. When $X$ is definably compact, $\Pow_n(X)$ is also definably compact.



Let $G$ be a definable group and $X$ a definable set. Recall that an action $\rho:G\times X\to X$ of $G$ on $X$ is definable over $k$ if $G$, $X$, and $\rho$ are all definable over $k$.

\begin{fact}[Lemma 2.23 of \cite{JY-1}]\label{fact-def-type-in-compact-space}
Let $ X $ be a definable topological space.   Let $ p $ be a definable $ 1 $-dimensional complete  type over $\M$ on $X$. If $X$ is definably compact, then $ p $ specializes to a point $c\in X $, meaning that $ p(x) \vdash x \in U $ for every  definable neighborhood $ U \ni c $.
\end{fact}

\begin{remark}
In Lemma 2.23 of \cite{JY-1}, $ X $ is a definable manifold, but in fact its proof applies to any definable compact space.
\end{remark}

\begin{proposition}\label{prop-dfg-act-on-dcp}
Let $(G,X)$ be a definable continuous action of a dfg group $G$ on a definable topological space $X$. If $X$ is definably compact, then there exists a point in $X$ with a finite orbit. Equivalently, there exists a natural number $n$ such that the induced action of $G$ on $\Pow_n(X)$ has a fixed point.
\end{proposition}

\begin{proof}
We proceed by induction on $\dim G$. For $g \in G$ and $x \in X$, let $g \cdot x$ denote the image of $x$ under the action of $g$.

Assume that $\dim(G) = 1$. Fix some $x_0\in X$. Consider the definable map $\tau \colon G \to X$ defined by $g \mapsto g \cdot x_0$. Then $G/\stab(x_0) = \tau(G)$ is definable, so $\dim(\tau(G)) = \dim(G) - \dim(\stab(x_0))$. If $\tau(G)$ is finite, then $x_0$ has a finite orbit, and we are done. Otherwise, $\stab(x_0)$ is finite, and hence $\tau$ is a finite-to-one map.

Let $p$ be a dfg type of $G$, and set $q = \tau(p)$. Then $q$ is a $1$-dimensional definable type on $X$. By Fact~\ref{fact-def-type-in-compact-space}, $q$ specializes to a point $c$ in $X$, meaning that $q \vdash (x \in U)$ for every definable open neighborhood $U$ of $c$. Since $X$ is Hausdorff, $c$ is unique. It is obvious that $g\cdot c$ is the unique point in $X$ that $g\cdot q = \tau(g\cdot p)$ specializes to. Since $\{g\cdot p \mid g\in G\}$ is bounded, $\{\tau(g\cdot p) \mid g\in G\}$ is bounded as $\tau$ is finite-to-one, and thus $\{g\cdot c \mid g\in G\}$ is also bounded. As the orbit $\{g\cdot c \mid g\in G\}$ is definable, it is finite.

Next, assume that $\dim(G) > 1$ and that the result holds for groups of lower dimension. By Theorem~\ref{thm-dfg-split}, we may assume $G$ has a normal definable dfg subgroup $H$ such that $G/H$ is a one-dimensional dfg group. By induction, there exists $n\in\mathbb{N}$ for which the set $Y$ of fixed points of $H$ in $\Pow_n(X)$ is non-empty. Clearly, $Y$ is definable and closed in $\Pow_n(X)$, and hence definably compact. Since $H$ is normal, $Y$ is closed under the $G$-action, so $G/H$ acts naturally on $Y$, and this action is definable. By induction, there exists $\ell\in\mathbb{N}$ such that the action of $G/H$ on $\Pow_\ell(Y)$ has a fixed point. It is then straightforward to verify that the original action of $G$ on $X$ has a point whose orbit size is at most $\ell n$. This completes the proof.
\end{proof}

\begin{theorem}[Conjugate-Commensurable Theorem]\label{theorem-VCT}
Let $G$ be a definable group, and let $H_1$ and $H_2$ be dfg subgroups such that $G/H_2$ is definably compact. Then there exists a finite index subgroup $A$ of $H_1$ contained in a conjugate of $H_2$. In particular, $H_2$ is a dfg component of $G$.
\end{theorem}
\begin{proof}
We let $H_1$ act naturally on $G/H_2$. By the preceding lemma, some $zH_2\in G/H_2$ has a finite orbit—namely, there exists a finite subset $Z$ of $G$ with $z\in Z$ such that $H_1 zH_2\subset Z H_2$. It follows that $H_1\subset Zz^{-1} \cdot z H_2 z^{-1}$. Let $W=Z z^{-1}$ and $H_3=z H_2 z^{-1}$; then $H_1\subset W H_3$. For each $w\in W$, if $H_1\cap w H_3\neq \emptyset$, we may replace $w$ with an element $w'\in H_1\cap w H_3$. Thus, we may assume $W\subset H_1$.

We then have
\[
H_1= \bigcup_{w\in W} (H_1\cap wH_3) =\bigcup_{w\in W} w \left(w^{-1}H_1\cap H_3\right)=\bigcup_{w\in W} w(H_1\cap H_3).
\]
Hence, $A =H_1\cap H_3$ is a finite index subgroup of $H_1$ contained in $H_3=z H_2 z^{-1}$.
\end{proof}

\section{Actions of $\dfg$ groups on $\fsg$ groups}\label{section-act-dfg-fsg}

In this section,  $H$ and $C$ will denote a dfg group and an $\fsg$ group defined in $\M$  with parameters from $k$, respectively, and  $\rho: H\times C\to C$ will be a $k$-definable action of $H$ on $C$, namely, $\rho(h,-):C\to C$ is an automorphsim of $C$ and $\rho(h,\rho(h',x))=\rho(h \cdot h',x) )$ for any $h,h'\in H$. Clear $\rho$ induces a  group morphism  $h\mapsto \rho(h,-)$ from $H$ to $\Aut(C)$, we denote it by $\rho^*: H\to \Aut(C)$. Clearly, each $\sigma\in\rho^*(H(k))$ is a $k$-definable automorphsim of $C$, so we may assume that $\rho^* (H(k))\leq \Aut(C(k))$.

\begin{lemma}\label{action-lemma}
If $k=\Qp$, then the image $\image \rho^*$ of $\rho^*$ is finite. 
\end{lemma}

\begin{proof}


Note that $C(\Q_p)$ is a compact $p$-adic analytic group (see p. 234 in \cite{OP}). 
By Corollary 8.35 in \cite{pro-p}, $\Aut(C(\Q_p))$ is a compact $p$-adic analytic group. Moreover, its topology is given by subgroups 
\[
\widetilde{N}:=\{\sigma\in \Aut(C(\Q_p)):g^{-1}\sigma(g)\in N \text{ for every } g\in C(\Q_p)\}
\]
as $N$ runs over the open normal subgroups of $C(\Q_p)$ (see p. 89 in \cite{pro-p}). Note that $C$ admits a $\Q_p$-definable family $\mathcal{U}=\{U_b: b\in Y\}$ of normal subgroups which forms a neighbourhood system of $\id_C$ (see Observation 8.6 in \cite{Johnson-P-minimal groups}). It is clear that the topology on $\Aut(C(\Q_p))$ can be also given by subgroups $\widetilde{U_b(\Q_p)}$ for $b\in Y(\Q_p)$. Note that those $\widetilde{U_b(\Q_p)}$ are open in $\Aut(C(\Q_p))$ and hence of finite index. So, ${\rho^*}^{-1}(\widetilde{U_b(\Q_p)})$ is a finite index subgroup of $H(\Q_p)$.

For $b\in Y$, we let $X_b:=\{h\in H: g^{-1}{\rho }(h,g)\subset U_b \text{ for every }g\in C\}$. 
Obviously, if $b\in Y(\Q_p)$, then $X_b(\Q_p)={\rho^*}^{-1}(\widetilde{U_b(\Q_p)})$ has finite index in $H(\Qp)$. Let $X=\bigcap_{y\in Y }X_y$, then $X$ is $\Qp$-definable. 
By Corollary \ref{coro-intersection-of-a-family-of-dfg}, we see that there are $b_1,...,b_m\in Y(\Qp)$ such that $X(\Qp)=\bigcap_{b\in Y(\Qp)}X_b(\Qp)=\bigcap_{i=1}^m X_{b_i}(\Qp)$. It follows that $X$ has finite index in $H$.


Note that $\ker {\rho^*}=\{h\in H: {\rho^*}(h,g)=g \text{ for every }g\in C\}$ is definable over $\Q_p$, and
\[
(\ker {\rho^*}) (\Q_p)=\ker ({\rho^*\upharpoonleft_{\Qp}})=\bigcap_{b\in Y(\Q_p)} {(\rho^*\upharpoonleft_{\Qp})}^{-1}(\widetilde{U_b(\Q_p)})=\bigcap_{b\in Y(\Q_p)} X_b(\Q_p)= X(\Q_p) .
\]
Hence, $\ker {\rho^*}= X$   is a finite index subgroup of $H$. Consequently, $\image {\rho^*}$ is finite.
\end{proof}

Suppose that $G$ is a semidirect product of $C$ and $H$. Then the conjugation
\[
\rho \colon H\times C\to C,\ (h,c)\mapsto c^h
\]
is a definable action of $H$ on $C$. Conversely, any definable action $\rho \colon H\times C\to C$ defines a semidirect product structure of $C$ and $H$. So we will use the notation $C\rtimes_{\rho} H$ to denote the semidirect product where the definable action of $H$ on $C$ is $\rho$. Clearly, $C\rtimes_{\rho} H$ is a direct product iff $\rho^* \colon H\to \operatorname{Aut}(C)$ is trivial.

As an immediate corollary of Lemma~\ref{action-lemma}, we have

\begin{corollary}\label{corollary-semi-product=product}
 If $k=\Qp$, then there is a finite-index subgroup $H_1$ of $H$ such that $C\rtimes_\rho H_1 = C\times H_1$.
\end{corollary}

We now show that the above corollary holds for any $k \models \pCF$.

\begin{remark}
Note that the following are equivalent:
\begin{enumerate}
    \item [(i)] there exists a finite-index subgroup $H_1 \leq H$ such that $C \rtimes_\rho H_1 = C \times H_1$;
    \item [(ii)] $H \cap Z(C)$ has finite index in $H$;
    \item [(iii)] $H \cap Z(C)$ is an open dfg subgroup of $H$,
\end{enumerate}
where $Z(C)$ denotes the centralizer of $C$, computed in $C \rtimes_\rho H$.
\end{remark}


\begin{proposition}\label{proposition-semi-product-of-fsg-dfg}
There is a $k$-definable finite-index subgroup $A \leq H$ such that $C \times A$ is a finite-index subgroup of $C \rtimes_\rho H$.
\end{proposition}

\begin{proof}
     Suppose that there is a counterexample in $k$, namely, there are $\emptyset$-definable families $\{C_y \mid y\in Y\}$, $\{H_y \mid y\in Y\}$, and $\{\rho_y \mid y\in Y\}$, and $a\in k^{|y|}$ such that $C_a$ has fsg, $H_a$ has dfg, and $\rho_a$ is a definable action of $H_a$ on $C_a$, and there is no finite-index subgroup $A$ of $H_a$ such that $C_a \rtimes_{\rho_a} A = C_a \times A$. By Fact~\ref{fact-def-cp-definable} and Theorem~\ref{thm-dfg-definable}, we may assume that $Y$ is $\emptyset$-definable. We see that
    \[
    k \models \exists y\in Y \left(Z(C_y) \cap H_y \text{ is not an open dfg subgroup of } H_y\right),
    \]
    which is an $\mathcal{L}$-sentence. So
    \[
    \Qp \models \exists y\in Y \left(Z(C_y) \cap H_y \text{ is not an open dfg subgroup of } H_y\right),
    \]
    and thus there is also a counterexample in $\Qp$. This is a contradiction.
\end{proof}

\section{The $\dfg$/$\fsg$ decompositions and $\dfg$ components}\label{section-dfg-fsg}



In this section, we assume that  $G$ is a group that is defined within the monster model $\M$  and is definable over the field $k$. We will show in this section  that $G$ has a $\dfg$/$\fsg$ decomposition when $G$ is definably amenable.

First, we prove the uniqueness of the $\dfg$/$\fsg$ decomposition up to commensurability, assuming such a decomposition exists:

\begin{lemma}\label{lemma-dfg-components-unique}
    Let $H_1$ and $H_2$ be two normal $\dfg$ subgroups of $G$ such that $G/H_1$ and $G/H_2$ are $\fsg$ groups. Then $H_1 \cap H_2$ has finite index in both $H_1$ and $H_2$.
\end{lemma}

\begin{proof}
    Note that $H_1 \cap H_2$ has dfg by Lemma~\ref{lemma-dfg-intersection}. By Fact~\ref{fact-open-dfg-subgroup}, both $H_2/(H_1 \cap H_2)$ and $G/H_1$ are definable groups. Since $H_2$ has dfg, $H_2/(H_1 \cap H_2)$ has dfg. Now $H_2/(H_1 \cap H_2) = H_2 H_1 / H_1$ is a definable subgroup of $G/H_1$, hence it is closed in $G/H_1$.

    As $G/H_1$ is definably compact, the closed subgroup $H_2/(H_1 \cap H_2)$ is also definably compact, and thus has fsg. Since it has both dfg and fsg, $H_2/(H_1 \cap H_2)$ must be finite. By symmetry, $H_1/(H_1 \cap H_2)$ is also finite.
\end{proof}

\begin{lemma}\label{find-normal}
Let $H\lhd U\lhd V$ be groups defined over $k$, where $H$ is a dfg group and $U/H$ is an fsg group. Then there exists a $k$-definable finite-index subgroup $\widetilde{H}$ of $H$ that is normal in $V$.
\end{lemma}

\begin{proof}
Pick an arbitrary $y\in V$. The conjugation $x\mapsto x^y$ is a definable automorphism of $U$, therefore $H^y\lhd U^y=U$ and $U/H^y$ has fsg.  It follows from   Lemma \ref{lemma-dfg-components-unique} that $H^y\cap H$ is a finite index subgroup of $H$ and equivalently, $H^0<H^y$.

Let $\widetilde{H}=\bigcap_{y\in V} H^y$. Obviously, $\widetilde{H}$ is a $k$-definable normal subgroup of $V$ containing $H^0$, which implies that $\widetilde{H}$ is a finite index subgroup of $H$.
\end{proof}

\subsection*{Existence in linear case}

We first study the case when $G$ is  linear, namely, that $G$ is a closed subgroup of $\GL_n(\M)$. 

\begin{notation}\label{notation}
Let $G$ be a linear group definable over $k$; let $\bar{G}$ denote the Zariski closure of $G$ (in $\Omega$); let $R$ denote the solvable radical of $\bar{G}$; let $S = \bar{G}/R$; and let $W = R(\mathbb{M}) \cap G$. Then $S$ is semi-simple, $W\lhd G$, and $G/W$ is open in $S(\mathbb{M})$. Note that $G$ is definably amenable if and only if both $W$ and $G/W$ are definably amenable (see Exercise 8.23 in \cite{Simon}).
\end{notation}

\begin{lemma}\label{lemma-G/W-cpt}
Let $G$ and $W$ be as in Notation \ref{notation}. If $G$ is definably amenable, then $G/W$ is definably compact.
\end{lemma}

\begin{proof}
By Fact \ref{fact-connected-semi-simples}, we can assume that $S=\bar G/R$ is $k$-simple.  Suppose for a contradiction that $G/W$ is not definably compact. Since $G/W$ is open in $S(\M)$, we see that   $S$ is not definably compact and hence is $k$-isotropic by Corollary \ref{coro-definably-compact=not-split}. By Corollary \ref{coro-G+=G0} (ii),   $G/W$ has finite index in $S(\M)$ and therefore $S(\M)$ is definably amenable, which contradicts Lemma \ref{lemma-simple+amenable is compact}.
\end{proof}

\begin{lemma}\label{lemma-W-definably-solvable}
Let $G$ and $W$ be as in Notation~\ref{notation}. Suppose that $G/W$ has fsg and $W$ has a $\dfg$/$\fsg$ decomposition. Then $G$ has a $\dfg$/$\fsg$ decomposition.
\end{lemma}

\begin{proof}
    Let $H$ be a normal dfg subgroup of $W$ such that $W/H$ has fsg. By Lemma~\ref{find-normal}, we may assume that $H\lhd G$. Now $G/H$ is an extension of $G/W$ by $W/H$, and both $G/W$ and $W/H$ have fsg. It follows from Remark~\ref{rmk-fsg-exact-seq} that $G/H$ also has fsg.
\end{proof}

Now we deal with $W$. It is clear that $W$ can be solved using the subnormal series of $R$ which is definable over $k$. Hence, we have a subnormal series  
\begin{equation}\label{equ-subnormal-series}
\{1_W\}=W_0\lhd W_1 \lhd\cdots\lhd W_n = W   
\end{equation}
defined over $k$ such that each $W_{i}/W_{i - 1}$ is $k$-definable and commutative. We call definable groups (need not be linear) with such definable series {\em definably solvable} (over $k$).

\begin{proposition}\label{proposition-def-solvable} 
Let $B$ be a definable group. If $B$ is definably solvable, then $B$ has a $\dfg$/$\fsg$ decomposition.
\end{proposition}
\begin{proof}
Let 
\[
\{1_B\}=B_0\lhd B_1 \lhd\cdots\lhd B_n = B
\]
be a subnormal series of $B$ such that each $B_{i}/B_{i - 1}$ is definable and commutative. For convenience, we assume $\dim (B_i/B_{i - 1})\geq 1$ for each $i$.

Recall that Johnson and Yao proved the dfg/fsg decomposition theorem for commutative definable groups in $\pCF$ (Theorem 1.1 in \cite{JY-2}).

We prove that $B$ has a dfg/fsg decomposition by induction on $\dim B$. If $\dim B=1$, the result is clear.

Now assume $\dim B>1$. If $n=1$, then $B=B_1$ is commutative, so $B$ has a dfg/fsg decomposition. Thus we assume $n>1$. By induction, $B_{n-1}$ and $B/B_{n-1}$ have $\dfg$/$\fsg$ decompositions:
\[
1 \rightarrow H_1 \rightarrow B_{n-1} \rightarrow C_1 \rightarrow 1,
\]
\[
1 \rightarrow H_2 \rightarrow B/B_{n-1} \rightarrow C_2 \rightarrow 1.
\]
By Lemma \ref{find-normal}, we may assume $H_1 \lhd B$. Let $\pi_1 : B\rightarrow B/H_1$ and $\pi_2 : B\rightarrow B/B_{n-1}$ be natural projections, and set $V:=\pi_1(\pi_2^{-1}(H_2))$. Clearly, $V$ is solvable. Since $H_1\leq \pi_2^{-1}(H_2)$, we have $\pi_1^{-1}(V)=\pi_2^{-1}(H_2)$.

\vspace{0.2cm}
\noindent \underline{Case 1: $\dim V<\dim B$:}

By induction, $V$ has a $\dfg$/$\fsg$ decomposition:
\[
1 \rightarrow H_3 \rightarrow V \rightarrow C_3 \rightarrow 1.
\]
By Lemma \ref{find-normal}, assume $H_3 \lhd B/H_1$. Let $H_4:=\pi_1^{-1}(H_3)$, so $H_4 \lhd B$. Since $H_4$ is an extension of dfg groups $H_1$ and $H_3$, it is dfg. Moreover,
\[
C_2\cong \frac{B/B_{n-1}}{H_2}\cong B/\pi_2^{-1}(H_2)\  =B/ \pi_1^{-1}(V)\cong \frac{B/H_4}{\pi_1^{-1}(V)/H_4}.
\]
Note that $\pi_1^{-1}(V)/H_4\cong V/H_3\cong C_3$. Thus $B/H_4$ is an extension of fsg groups $C_3$ and $C_2$, hence fsg. Thus $B$ admits the $\dfg$/$\fsg$ decomposition
\[
1 \rightarrow H_4 \rightarrow B \rightarrow B/H_4 \rightarrow 1.
\]

\noindent \underline{Case 2: $\dim V =\dim B$:}

Here, $\dim H_1=\dim C_2=0$. For convenience, assume $B_{n-1}=C_1$ and $B/B_{n-1}=B/C_1=H_2$, so $\pi_2: B\rightarrow H_2$ is surjective with $\ker \pi_2=C_1$.

If $B$ has fsg, we are done. Otherwise, by the Peterzil–Steinhorn Theorem in $\pCF$ (Corollary 6.11 in \cite{JY-3}), there is a one-dimensional dfg subgroup $U$ of $B$. Now $C_1U$ is a definable subgroup of $B$. Consider the action of $U$ on $C_1$ by conjugation. By Lemma \ref{action-lemma}, replacing $U$ with a finite-index subgroup, we may assume $U \leq Z_{C_1U}(C_1)$. Thus $U \lhd C_1U$, and $(C_1U)/U$ is definable with fsg (since $C_1$ has fsg). Hence
\[
1 \rightarrow U \rightarrow C_1 U \rightarrow (C_1U)/U \rightarrow 1
\]
is a decomposition of $C_1U$.

Since $H_2=B/C_1$ is commutative, $\pi_2(U) \lhd H_2$, so $C_1 U=\pi_2^{-1}(\pi_2(U)) \lhd B$. Moreover, $\pi_2(U)=(UC_1)/C_1$ is a dfg subgroup of $H_2$, so   $B/(C_1U)\cong H_2/\pi_2(U)$ is definable and commutative. Replacing $B_{n-1}$ with $C_1U$ reduces to Case 1.
\end{proof}

We conclude from Lemma\ref{lemma-G/W-cpt}, Lemma \ref{lemma-W-definably-solvable} and Proposition \ref{proposition-def-solvable} that 

\begin{corollary}\label{cor-linear and definably amenable}
   Let $G$ and $W$ be as in Notation \ref{notation}, then the following are equivalent
    \begin{enumerate}
        \item [(i)] $G$ is definably amenable;
        \item [(ii)] $G/W$ is definably compact;
        \item [(iii)] $G$ has a $\dfg$/$\fsg$ decomposition.
    \end{enumerate} 
\end{corollary}

\subsection*{Existence in general case}

We now present the general version.

\begin{theorem}\label{theorem-DA-has-dfg-fsg-decomp}
Let $G$ be an arbitrary definable group over $k$. If $G$ is definably amenable, then $G$ has a $\dfg$/$\fsg$ decomposition.
\end{theorem}
\begin{proof}
By Fact \ref{fact-antilinear-linear}, there is a $k$-definable exact sequence 
\[
1\rightarrow A\rightarrow G \stackrel{\pi_G}{\rightarrow} V\rightarrow 1,
\]
where $A$ is commutative-by-finite and $V$ is linear. The center $C(A)$ of $A$ is also a definable normal subgroup of $G$. Replacing $A$ by $C(A)$, we may assume that $A$ is commutative and $V$ contains a finite subgroup $X$ such that $\bar V=V/X$ is a definably amenable linear group. Thus by Corollary \ref{cor-linear and definably amenable}, $\bar V$ has a normal dfg subgroup $\bar H$ such that $\bar V/\bar H$ is an fsg group. Let $\pi_V: V\to \bar V$ be the natural projection. Then $H=\pi_V^{-1}(\bar H)$ is dfg and $C=V/H\cong \bar V/\bar H$ is fsg, so $V$ has a $\dfg$/$\fsg$ decomposition 
\[
1\to H\to V\to C\to 1.
\]
The subgroup $\pi_G^{-1}(H)$ is normal in $G$. Moreover, $\pi_G^{-1}(H)$ is definably solvable as it is an extension of $H$ by $A$, and both $H$ and $A$ are definably solvable. By Proposition \ref{proposition-def-solvable}, $\pi_G^{-1}(H)$ has a $\dfg$/$\fsg$ decomposition
\[
1 \rightarrow H' \rightarrow \pi_G^{-1}(H) \rightarrow C'\rightarrow 1. 
\]
Note that
\[
C\cong V/H\cong G/\pi_G^{-1}(H)\cong \frac{G/H'}{\pi_G^{-1}(H)/H'},
\]
where $\pi_G^{-1}(H)/H'\cong C'$. It follows that $G/H'$ is an extension of the fsg groups $C'$ and $C$, and hence $G/H'$ is an fsg group. We conclude that 
\[
1 \rightarrow H' \rightarrow G \rightarrow G/H'\rightarrow 1
\]
is a $\dfg$/$\fsg$ decomposition of $G$.
\end{proof}

\begin{corollary}\label{coro-dfg-component-of-DA}
Let $G$ be an arbitrary definable group over $k$. If $G$ is definably amenable, and $H_1, H_2$ are dfg components of $G$, then $H_1 \cap H_2$ has finite index in both $H_1$ and $H_2$. In particular, $\bigcap_{g\in G} H_1^g$ is a normal dfg component of $G$.
\end{corollary}

\begin{proof}
By Theorem \ref{theorem-DA-has-dfg-fsg-decomp}, $G$ has a normal dfg component $A$ such that $G/A$ is definably compact. Since $H_1/A = H_1A/A$ is a closed subgroup of $G/A$, it is also definably compact. On the other hand, $H_1/A$ has dfg because $H_1$ has dfg. We conclude that $H_1/A$ has both dfg and fsg, so it is finite. Similarly, $H_2/A$ is finite. Thus, both $H_1 \cap A$ and $H_2 \cap A$ are open dfg subgroups of $A$ (as they have the same dimension), so by Fact \ref{fact-open-dfg-subgroup}, $H_1 \cap A$ and $H_2 \cap A$ have finite index in $A$. It follows that $H_1 \cap H_2 \cap A$ has finite index in $A$. Hence, $H_1 \cap H_2$ is an open dfg subgroup of both $H_1$ and $H_2$, and thus has finite index in both.

It is easy to see that $\bigcap_{g\in G} H_1^g$ is normal. Since $H^0 \leq \bigcap_{g\in G} H_1^g$, we see that $\bigcap_{g\in G} H_1^g$ has finite index in $H_1$, so it is also a dfg component of $G$.
\end{proof}

\begin{remark}
 Recall that we work in the monster model $\M$. As defined in the introduction, a dfg component of a definable group $G$ is a dfg subgroup of $G$ of maximal dimension. We adopt the monster model as the ambient setting for this definition to ensure it is well-posed: in $\M$, a maximal-dimensional dfg subgroup always exists.
However, by the definability of both the dfg property (Theorem~\ref{thm-dfg-definable}) and dimension, the notion of a dfg component is first-order characterizable. It follows that whenever $G$ is definable over $k$, $G$ admits a dfg component that is itself definable over $k$ (see Corollary \ref{coro-dfg component over k}).
\end{remark}
\subsection*{The $\dfg$ components and quotients}



Let us now adopt the perspective of the quotient $G/H$. Recall that by Fact \ref{fact-open-dfg-subgroup}, $G/H$ is definable when $H$ is a dfg subgroup of $G$. We now give more results on the quotient spaces.

\begin{lemma}\label{lemma-dfg-components-over Qp}
Suppose that $ E \leq D \leq G $ are $\Qp$-definable groups such that $ G/E $, $ D/E $, and $ G/D $ are definable. If $ D/E $ and $ G/D $ are definably compact, then $ G/E $ is definably compact. 
\end{lemma}
\begin{proof}
 Let $ C(\Qp) $ be an open compact subgroup of $ G(\Qp) $. Let $ \pi \colon G \to G/D $ be the natural projection; since $\pi$ is an open map, $\pi(C(\Qp))$ is open in $G(\Qp)/D(\Qp)$. As $ G(\Qp)/D(\Qp) $ is compact, there exist $ g_1, \dots, g_n \in G(\Qp) $ such that  
\[ G(\Qp)/D(\Qp) = g_1\pi(C(\Qp)) \cup \dots \cup g_n\pi(C(\Qp)), \]  
which implies $ G(\Qp) = g_1CD \cup \dots \cup g_nCD $. Similarly, by the compactness of $ D(\Qp)/E(\Qp) $, there exist $ h_1, \dots, h_m \in D(\Qp) $ such that  
\[ D(\Qp) = h_1(C \cap D)E \cup \dots \cup h_m(C \cap D)E. \]  
Combining these, we get  $G(\Qp) = \bigcup_{i=1}^n \bigcup_{j=1}^m g_i h_j CE$.  
It follows that $ G/E = \bigcup_{i=1}^n \bigcup_{j=1}^m g_i h_j (C/E) $, which is definably compact because $ C/E $ is definably compact (as it is the image of a definably compact set).
\end{proof}

We now show that Lemma \ref{lemma-dfg-components-over Qp} is true for any $k$.

\begin{lemma}\label{lemma-dfg-components-over k}
Retain the notations $G$, $D$, and $E$ introduced in Lemma \ref{lemma-dfg-components-over Qp}, with $\Qp$ replaced by an arbitrary $p$-adically closed field $k$. Then 
\begin{enumerate}
    \item [(i)] If $D/E$ and $G/D$ are definably compact, then $G/E$ is definably compact.
    \item [(ii)] 
    If $ G/E $ is definably compact, then $ D/E $ is definably compact.
\end{enumerate} 
\end{lemma}
\begin{proof}
(i) Assume there exist definable families of groups $\{G_b\mid b\in B\}$, $\{D_b\mid b\in B\}$, $\{E_b\mid b\in B\}$ such that:
\begin{itemize}
    \item $E_b\leq D_b\leq G_b$ for each $b\in B$;
    \item $G=G_a$, $D=D_a$, $E=E_a$ for some $a\in B$.
\end{itemize}
Since $G_a/D_a$ and $D_a/E_a$ are definable and definably compact, and $G_a/E_a$ is definable, shrink $B$ if necessary to assume there exist:
\begin{itemize}
\item  Definable families of spaces $\{X_b\mid b\in B\}$, $\{Y_b\mid b\in B\}$, $\{Z_b\mid b\in B\}$;
\item  Definable families of continuous open surjective maps $\{f_b: G_b \to X_b\mid b\in B\}$, $\{g_b: D_b \to Y_b\mid b\in B\}$, $\{h_b: G_b \to Z_b\mid b\in B\}$;
\end{itemize}

satisfying for each $b\in B$:
\begin{itemize}
    \item $X_b$ and $Y_b$ are definably compact;  
    \item $f_b^{-1}(x)$ is a coset of $D_b$ in $G_b$ for all $x\in X_b$;
    \item $g_b^{-1}(y)$ is a coset of $E_b$ in $D_b$ for all $y\in Y_b$;
    \item $h_b^{-1}(z)$ is a coset of $E_b$ in $G_b$ for all $z\in Z_b$.
\end{itemize}
By Lemma \ref{lemma-dfg-components-over Qp}, we have:
\[
\Qp\models \forall b\in B(\Qp)\ \bigl(Z_b(\Qp)\text{ is definably compact}\bigr),
\]
which transfers to $\M$. Thus $G_a/E_a$ is definably compact as required.

(ii) Since $\pi \colon G\to G/E$ is an open map, $D$ is closed in $G$, and $E\leq D$, we have that $D/E = \pi(D)$ is closed in $G/E$. Thus $D/E$ is definably compact, as a closed subset of the definably compact space $G/E$.
\end{proof}

    
\begin{lemma}\label{lemma-G/dfg-component-def-cp}
    If $G$ has a $\dfg$ subgroup $H$ such that $G/H$ is definably compact. Then for any $\dfg$ component $H'$ of $G$, $G/H'$ is definably compact.
\end{lemma}
\begin{proof}
 By Theorem \ref{theorem-VCT}, there is a finite index subgroup $H_1$ of $H' $  and some $g\in G$ such that $H_1^g\leq H$.  Since $H_1$ is also a dfg component of $G$, we have 
 \[
 \dim(H)=\dim(H_1^g)=\dim(H_1)=\dim(H'),\]
 which means that $H/H_1^g$ is finite and thus is definably compact.  
 It follows from Lemma \ref{lemma-dfg-components-over k} that $G/H_1^g$ is definably compact. So $G/H_1$ is definably compact, and therefore $G/H'$ is also definably compact.   
\end{proof}

For our purposes, we generalize Fact~\ref{fact-quotient-by-max-split-subgroup} to any $p$-adically closed field  $k$ for semi-simple algebraic groups.

\begin{fact}[\cite{ABC}]\label{fact-maximal solvable $k$–subgroups-conjugation}
    Let $ G $ be a semi-simple algebraic group over $ k $. Then all maximal solvable $ k $-subgroups of $ G $ are conjugate by elements of $ G(k) $.
\end{fact}

\begin{lemma}\label{lemma-split-solvable-compact}
    Let $ P \leq G $ be linear algebraic groups over $ k $. If $ G $ is centreless and semi-simple, then $ G(k)/P(k) $ is definably compact if and only if $ P $ contains a maximal $ k $-split solvable subgroup of $ G $.
\end{lemma}
\begin{proof}
    By Theorem~\ref{Thm-semi-simple and K-simple}, $ G $ is $ k $-isomorphic to a semi-simple algebraic group $ G' $ over $ \mathbb{Q}_p $ via $ f \colon G \to G' $. Let $ H' \leq G' $ be a maximal $ \mathbb{Q}_p $-split solvable subgroup of $ G' $; then $ H = f^{-1}(H') $ is a $ k $-split solvable subgroup of $ G $ (as $ H' $ $ k $-embeds into a trigonalizable algebraic group over $ k $). By Fact~\ref{fact-quotient-by-max-split-subgroup}, $ G'(k)/H'(k) $ is definably compact, so $ G(k)/H(k) $ is also definably compact. Thus $ H(k) $ is a dfg component of $ G(k) $ by Theorem~\ref{theorem-VCT}, and hence $ H $ is a maximal $ k $-split solvable subgroup of $ G $ by Lemma~\ref{lemma-linear-alg-group-k-split-iff-dfg}. 
    
    Suppose that $ P \leq G $ is an algebraic group over $ k $ containing a maximal $ k $-split solvable subgroup $ A $. By Fact~\ref{fact-maximal solvable $k$–subgroups-conjugation}, $ H $ is conjugate to $ A $ by an element of $ G(k) $, so $ G(k)/A(k) $ is definably compact. The natural projection $ G(k)/A(k) \to G(k)/P(k) $ then implies that $ G(k)/P(k) $ is also definably compact.
    
    Conversely, suppose $ G(k)/P(k) $ is definably compact. Let $ H $ be a $ k $-split solvable subgroup of $ G $ such that $ G(k)/H(k) $ is definably compact. Let $ G(x,y) $, $ P(x,y) $, $ H(x,y) $ be quantifier-free formulas and $ a \in k^{|y|} $ such that $ G = G(\Omega,a) \leq \GL_n(\Omega) $, $ P = P(\Omega,a) $, $ H = H(\Omega,a) $. By Fact~\ref{fact-quotient-by-max-split-subgroup},
    \[
    \begin{split}
    \mathbb{Q}_p \models \forall y \bigg( &\big(G(\Qp,y)\text{ is a subgroup of } \GL_n(\Qp)\big)\wedge \big(H(\mathbb{Q}_p,y) \text{ is a dfg subgroup of } G(\mathbb{Q}_p,y)\big) \wedge \\
    &\big(P(\mathbb{Q}_p,y) \text{ is a subgroup of } G(\mathbb{Q}_p,y)\big) \wedge \big(G(\mathbb{Q}_p,y)/H(\mathbb{Q}_p,y) \text{ is definably compact}\big)\wedge \\
    &\big(G(\mathbb{Q}_p,y)/P(\mathbb{Q}_p,y) \text{ is definably compact}\big)\bigg) \to \bigg(\exists g \in G(\mathbb{Q}_p) \big(H(\mathbb{Q}_p,y)^g \leq P(\mathbb{Q}_p,y)\big)\bigg),
    \end{split}
    \]
    which is a first-order sentence and hence transfers to $ k $. Thus there exists $ g \in G(k) $ such that $ H(k)^g \leq P(k) $, so $ H^g \leq P $; that is, $ P $ contains a maximal $ k $-split solvable subgroup of $ G $.
\end{proof}

Now we characterize dfg components using the quotient $G/H$.

\begin{lemma}\label{lemma-dfg-components-linear}
    Let $G$ be a definable group. If $ G $ is linear and $ H $ is a dfg component of $ G $, then $ G/H $ is definably compact.
\end{lemma}
\begin{proof}
Let $ \bar{G} $ be the Zariski closure of $ G=G(\M) $ in $ \GL_n(\Omega) $, $ R $ the solvable radical of $ \bar{G} $, $ S=\bar{G}/R $, and $ W=R(\M) \cap G $. Let $ S_1=G/W $ and $ \pi \colon G\to S_1 $ be the natural projection. Then $ S_1 $ is open in $ S(\M) $.

We distinguish two cases. If $S$ is $k$-anisotropic, then by Corollary~\ref{coro-definably-compact=not-split}, $ S(\M) $ is definably compact, and thus $ S_1 $ is definably compact as it is closed in $ S(\M) $. If $S$ is $k$-isotropic, then by Corollary~\ref{coro-G+=G0}, $ S_1 $ is either definably compact or of finite index in $ S(\M) $. In summary, $ S_1 $ is either definably compact or of finite index in $ S(\M) $.

If $ S_1 $ is definably compact, then the short exact sequence $ 1\to W\to G\to S_1\to 1 $ shows that $ G $ is definably amenable. By Theorem~\ref{theorem-DA-has-dfg-fsg-decomp}, $G$ admits a normal dfg subgroup $H'$ such that $G/H'$ is definably compact. By Lemma~\ref{lemma-G/dfg-component-def-cp}, $ G/H $ is also definably compact.

If $ S_1 $ is of finite index in $ S $, we let $ H_1 $ be a maximal $ k $-split solvable algebraic subgroup of $ S $. Then $ E=S_1\cap H_1(\M) $ has finite index in $ H_1(\M) $, so $ E $ is a dfg subgroup of $ S_1 $. Since both $ S(\M)/H_1(\M) $ and $ H_1(\M)/E $ are definably compact, by Lemma~\ref{lemma-dfg-components-over k}(i), $ S(\M)/E $ is definably compact. By Lemma~\ref{lemma-dfg-components-over k}(ii), $ S_1/E $ is also definably compact.

Let $ G_1=\pi^{-1}(E) $. Then $ G_1 $ is an extension of $ E $ by $ W $, hence definably amenable. By Theorem~\ref{theorem-DA-has-dfg-fsg-decomp}, $ G_1 $ has a $\dfg$/$\fsg$ decomposition $ 1\to H_2\to G_1\to C\to 1 $. Since both $ G/G_1\cong S_1/E $ and $ G_1/H_2 $ are definably compact, Lemma~\ref{lemma-dfg-components-over k}(i) implies $ G/H_2 $ is definably compact. By Lemma~\ref{lemma-G/dfg-component-def-cp}, $ G/H $ is also definably compact.
\end{proof}

\begin{proposition}\label{proposition-G/dfg-component-def-cp}
   For any definable group $ G $, if $ H $ is a dfg component of $ G $, then $ G/H $ is definably compact. 
\end{proposition}
\begin{proof}
There is a $ k $-definable exact sequence 
\[
1\rightarrow A\rightarrow G \stackrel{\pi}{\rightarrow} V\rightarrow 1,
\]
where $ A $ is commutative-by-finite and $ V $ is linear. By Lemma \ref{lemma-dfg-components-linear}, $ V/H_0 $ is definably compact for any dfg component $ H_0 $ of $ V $. Let $ D=\pi^{-1}(H_0) $; then $ G/D $ is definably homeomorphic to $ V/H_0 $, hence definable and definably compact.  

Since $ D $ is an extension of definably amenable groups, Theorem \ref{theorem-DA-has-dfg-fsg-decomp} implies $ D $ has a normal dfg component $ H_1 $. By Lemma \ref{lemma-dfg-components-over k} (i), $ G/H_1 $ is definably compact. Lemma \ref{lemma-G/dfg-component-def-cp} then gives $ G/H $ is definably compact.
\end{proof}

We conclude from Theorem~\ref{theorem-VCT} and Proposition~\ref{proposition-G/dfg-component-def-cp} that

\begin{theorem}\label{theorem-G/dfg-component=quotient def-cp}
Let $G$ be a definable group and $H$ a dfg subgroup of $G$. Then $H$ is a dfg component of $G$ iff $G/H$ is definably compact.   
\end{theorem}

Replacing a dfg component by its definable connected component, we can drop the ``commensurable'' part from the Conjugate-Commensurable Theorem (Theorem~\ref{theorem-VCT}):
\begin{corollary}\label{coro-Conjugate-Theorem}
Let $G$ be a definable group, and let $H_1, H_2$ be dfg components of $G$. Then there exists $g\in G$ such that $(H_1^0)^g = H_2^0$.
\end{corollary}
\begin{proof}
By Theorem~\ref{theorem-VCT} and Theorem~\ref{theorem-G/dfg-component=quotient def-cp}, there exists $g\in G$ such that $H_1^g \cap H_2$ has finite index in $H_2$. Since $\dim(H_1^g \cap H_2) = \dim(H_2) = \dim(H_1^g)$, the intersection $H_1^g \cap H_2$ is an open dfg subgroup of $H_1^g$. By Fact~\ref{fact-open-dfg-subgroup}, it also has finite index in $H_1^g$.

It follows that $(H_1^g)^0 = (H_1^g \cap H_2)^0 = H_2^0$. Since conjugation by $g$ is a definable group isomorphism, we have $(H_1^g)^0 = (H_1^0)^g$. Thus $(H_1^0)^g = H_2^0$, as required.
\end{proof}

\begin{corollary}\label{coro-dfg component over k}
  Let $G$ be a definable group. If $G$ is definable over $k$, then $G$ has a dfg component definable over $k$. 
\end{corollary}
\begin{proof}
Take any dfg component $H$ of $G$. Suppose that $\{H_b \mid b\in B\}$ is a definable family of dfg groups and $H=H_a$ for some $a\in B$. Then $\M\models \exists y\in B\,(G/H_y\ \text{is definably compact})$, so $k\models \exists y\in B(k)\,(G/H_y\ \text{is definably compact})$. It follows that $G$ has a dfg component definable over $k$.
\end{proof}

The following corollary extends Fact~\ref{fact-quotient-by-max-split-subgroup} from linear algebraic groups over $\Qp$ to definable groups over an arbitrary $p$-adically closed field $k$:
\begin{corollary}\label{corollary-G/E is compact iff E contains dfg cp}
  Let $G$ be a definable group, $E$ a definable subgroup of $G$ such that $G/E$ is definable. Then $G/E$ is definably compact iff $E$ contains a dfg component of $G$.  
\end{corollary}
\begin{proof}
Suppose that $G/E$ is definably compact. Let $H$ be a dfg component of $E$. Then $E/H$ is definably compact. Thus $G/H$ is definably compact by Lemma~\ref{lemma-dfg-components-over k}. It follows from Theorem~\ref{theorem-G/dfg-component=quotient def-cp} that $H$ is a dfg component of $G$.

Conversely, if $H\leq E$ is a dfg component of $G$, then $G/H$ is definably compact, and hence so is $G/E$.
\end{proof}

Combining Lemma~\ref{lemma-linear-alg-group-k-split-iff-dfg} and Corollary~\ref{corollary-G/E is compact iff E contains dfg cp}, for algebraic groups, we may extend Fact~\ref{fact-quotient-by-max-split-subgroup} from $ \mathbb{Q}_p $ to an arbitrary $ p $-adically closed field $ k $ as follows:

\begin{corollary}
    Let $ G $ be an algebraic group defined over $ k $, and $ E $ an algebraic subgroup of $ G $ defined over $ k $. Then $ G(k)/E(k) $ is definably compact if and only if $ E $ contains a maximal $ k $-split solvable algebraic subgroup of $ G $.
\end{corollary}

As another application of Theorem \ref{theorem-G/dfg-component=quotient def-cp}, we can characterize dfg groups by the property described in Proposition \ref{prop-dfg-act-on-dcp}, which shows that dfg groups are precisely those groups that have no interaction with definably compact spaces.

\begin{corollary}
Let $G$ be a definable group. Then $G$ is a dfg group iff every definable continuous action of $G$ on  a definably compact space $X$ admits a finite orbit.
\end{corollary}
\begin{proof}
One direction is Proposition \ref{prop-dfg-act-on-dcp}, so it suffices to prove the converse. Let $H$ be a dfg component of $G$ and $X=G/H$. By Theorem \ref{theorem-G/dfg-component=quotient def-cp}, $X$ is definably compact. Since the action is transitive and admits a finite orbit, $X$ is finite. Then $G$ admits a finite index dfg subgroup $H$, so $G$ is a dfg group.
\end{proof}

\section{The definably amenable components}\label{section-DAC}



We conclude our paper with the study of definably amenable components. Recall that a \emph{definably amenable component} of a definable group $G$ is a maximal definably amenable subgroup containing a dfg component of $G$.


\begin{lemma}\label{lemma-91}
   Let $G$ be a definable group and let $H$ be a dfg subgroup of $G$. Then $N_G(H^0)=\{g\in G\mid (H^0)^g=H^0\}$ is definable.    
\end{lemma}
\begin{proof}
  Let $g\in G$. As shown in the proof of Corollary~\ref{coro-Conjugate-Theorem}, $(H^0)^g=H^0$ if and only if $H\cap H^g$ has finite index in both $H$ and $H^g$. By Lemma~\ref{lemma-dfg-intersection}, $H\cap H^g$ is dfg. By Fact~\ref{fact-open-dfg-subgroup}, the dfg group $H\cap H^g$ has finite index in both $H$ and $H^g$ if and only if it is open in $H$, which is a first-order condition. Thus $N_G(H^0)=\{g\in G\mid H\cap H^g \text{ is an open subgroup of } H\}$, which is definable.
\end{proof}


\begin{lemma}\label{lemma-bigcap-conjugate}
Let $G$ be a definable group,  $H $   a dfg subgroup of $G$, and $H_1=\bigcap_{g\in N_G(H^0)}H^g$, then $N_G(H^0)=N_G(H_1)$.
\end{lemma}
\begin{proof}
    Clearly, $H_1^0=H^0$.  If $g\in N_G(H_1)$, then $h\mapsto h^g$ is a definable automorphism of $H_1$, so it fixes the $H_1^0=H^0$, so $g\in N_G(H^0)$. Conversely, if $g\in N_G(H^0)$, then ${H_1}^g=(\bigcap_{b\in N_G(H^0)}H^b)^g=\bigcap_{b\in N_G(H^0)}H^{bg}={H_1}$, so $g\in N_G(H_1)$.
\end{proof}

The following theorem indicates that $N_G(H^0)$ is a definably amenable component of $ G $ when $H$ is a dfg component of $G$.

\begin{theorem}\label{theorem-max-DA-subgroup}
Let $ G $ be a definable group, $ H $ a dfg component of $ G $, and $ B = N_G(H^0) $. Then:
\begin{enumerate}
    \item [(i)] $ B $ is definably amenable and $ B = N_G(B^{00}) $. 
    \item [(ii)] $B$ is a definably amenable component. In fact, any definably amenable component of $G$ is of the form $N_G(H_1^0)$ where $H_1$ is a dfg component of $G$.
    \item [(iii)] Any two definably amenable components of $ G $ are conjugate.
    \item [(iv)] $ G $ is definably amenable if and only if $ G = B $.
\end{enumerate}
\end{theorem}

\begin{proof}
(i) By Lemma~\ref{lemma-bigcap-conjugate}, $ B = N_G(H_1) $, where $ H_1 = \bigcap_{g \in N_G(H^0)} H^g $ is also a dfg component of $ G $. Note that $ B/H_1 $ is definably compact, as it is a closed subset of $ G/H_1 $ and $ G/H_1 $ is definably compact. Thus $ B $ is definably amenable.

To show $ B = N_G(B^{00}) $: clearly, $ B \leq N_G(B) \leq N_G(B^{00}) $, so it suffices to prove $ N_G(B^{00}) \leq B $. Let $ g \in N_G(B^{00}) $; then $ (H^0)^g=(H^{00})^g \leq (B^{00})^g = B^{00} $. Since $ H $ is a dfg component of $ B $, Corollary~\ref{coro-dfg-component-of-DA} implies $ H \cap H^g $ has finite index in both $ H $ and $ H^g $. Hence $ (H^0)^g = H^0 $, so $ g \in N_G(H^0) = B $.

(ii) Suppose that $B_1$ is definably amenable and contains $B$. Then $ B_1 $ has a normal dfg component $ H_1 $. Now $ H $ is also a dfg component of $ B_1 $; we see from Corollary~\ref{coro-dfg-component-of-DA} that $ H^0 = H_1^0 $. Thus $ B_1 \leq N_G(H_1) \leq N_G(H^0) = B $.

To prove the ``In fact'' claim, suppose $ B_1 $ is a definably amenable component, and let $H_1$ be a normal dfg component of $B_1$. Then $H_1$ is also a dfg component of $G$. By (i), $N_G(H_1^0)$ is definably amenable and contains $B_1$, so $B_1 = N_G(H_1^0)$.

(iii) Let $B_1, B_2$ be definably amenable components of $G$. By (ii), $B_1 = N_G(H_1^0)$ and $B_2 = N_G(H_2^0)$ for some dfg components $H_1$ and $H_2$ of $G$. By Corollary~\ref{coro-Conjugate-Theorem}, there exists $g\in G$ such that $(H_1^0)^{g^{-1}} = H_2^0$. Thus
\[
B_1^g = N_G(H_1^0)^g = N_G\!\left((H_1^0)^{g^{-1}}\right) = N_G(H_2^0) = B_2.
\]

(iv) Suppose $ G $ is definably amenable. By Lemma~\ref{lemma-bigcap-conjugate}, $ H^0 $ is normal in $ G $, so $ G = N_G(H^0) = B $.
\end{proof}
 
\begin{corollary}
  Let $G$ be a definable group. If $G$ is definable over $k$, then $G$ has a definably amenable component definable over $k$. 
\end{corollary}
\begin{proof}
    By Corollary~\ref{coro-dfg component over k}, $G$ has a dfg component $H$ definable over $k$. By Theorem~\ref{theorem-max-DA-subgroup}, $B = N_G(H^0)$ is a definably amenable component of $G$ definable over $k$.
\end{proof}

\begin{remark}\label{rmk-dim-DA-Not max}
While each definably amenable component of $ G $ is maximal among definably amenable subgroups containing a dfg component of $ G $, it need not be the largest in dimension: when $ G $ is definable over $ \Qp $, it contains a definably compact open subgroup $ O $, which is definably amenable and satisfies $ \dim(O)=\dim(G) $.     
\end{remark}

In the case where a definably amenable component is fully dimensional, we have:
\begin{corollary}
    Let $ G $ be a definable group and $ B $ a definably amenable component of $ G $, then $ G=B $ if and only if $ \dim(G)=\dim(B) $.
\end{corollary}
\begin{proof}
One direction is obvious, so it suffices to show the other. 

Suppose $ \dim(G)=\dim(B) $. Since $ B $ is a definable subgroup of $ G $, $ \dim(B)=\dim(G) $ implies $ B $ is open in $ G $. Let $ H $ be a dfg component of $ G $.  Without loss of generality, assume $ N_G(H)=N_G(H^0)=B $. And, to obtain a contradiction, we assume that $G\neq B$.

Let $\varphi(x;y),\psi(x;y)\in \mathcal{L}$ and $b\in\M$ such that $G$ is defined by $\varphi(x;b)$ and $H$ is defined by $\psi(x;b)$. 
Let $Y$ be the set of those $b'$ such that 
\begin{itemize}
    \item $H_{b'}:=\psi(x;b')$ is a $\dfg$ subgroup of the definable group $G_{b'}:=\varphi(x;b')$;  
    \item Let $B_{b'}:=N_{G_{b'}}(H_{b'})$. $B_{b'}=N_{G_{b'}}(H_{b'}^0)$;
    \item $ \dim(B_{b'})=\dim(G_{b'}) $;
    \item $G_{b'}\neq B_{b'}$.
\end{itemize}
Then by Theorem \ref{thm-dfg-definable}, Lemma \ref{lemma-91} and the definability of dimension, $Y$ is $\emptyset$-definable. Obviously, $b\in Y$, so $Y$ is nonempty in every model. 
Hence we can find a counterexample in $\Q_p$. We may thus assume $b\in Y(\Q_p)$ and $G$, $H$, $B$ are definable over $\Q_p$.

    As the quotient map $G\rightarrow G/H$ is open, $ B/H $ is open in $ G/H $. Since $ G(\Qp)/H(\Qp) $ is compact, finitely many $ G(\Qp) $-translates of $ B(\Qp)/H(\Qp) $ cover $ G(\Qp)/H(\Qp) $. This implies finitely many $ G(\Qp) $-translates of $ B(\Qp) $ cover $ G(\Qp) $, so $ G/B $ is finite, so $ G $ is definably amenable. By Theorem \ref{theorem-max-DA-subgroup} (iv), $ G=B $, which is a contradiction. 
\end{proof}

\end{document}